\input amstex\documentstyle{amsppt}  
\pagewidth{12.5cm}\pageheight{19cm}\magnification\magstep1
\NoBlackBoxes
\topmatter
\title Unipotent almost characters of simple $p$-adic groups\endtitle
\author G. Lusztig\endauthor
\address{Department of Mathematics, M.I.T., Cambridge, MA 02139}\endaddress
\dedicatory{To G\'erard Laumon, on his 60th birthday.}\enddedicatory
\thanks{Supported in part by National Science Foundation grant DMS-0758262.}\endthanks
\endtopmatter   
\document

\define\Irr{\text{\rm Irr}}
\define\sneq{\subsetneqq}

\define\tcs{\ti{\cs}}

\define\uE{\un E}

\define\da{\dagger}

\define\frl{\forall}

\define\si{\sim}

\define\sqc{\sqcup}

\define\qua{\quad}

\define\hE{\hat E}

\define\bY{\bar Y}

\define\lb{\linebreak}

\define\op{\oplus}

\define\part{\partial}
\define\em{\emptyset}

\define\iy{\infty}
\define\m{\mapsto}
\define\do{\dots}

\define\lra{\leftrightarrow}

\define\sm{\smallmatrix}
\define\esm{\endsmallmatrix}
\define\sub{\subset}    

\define\T{\times}
\define\ti{\tilde}
\define\nl{\newline}
\redefine\i{^{-1}}

\define\un{\underline}
\define\ov{\overline}
\define\ot{\otimes}
\define\bbq{\bar{\QQ}_l}

\define\Ad{\text{\rm Ad}}
\define\Hom{\text{\rm Hom}}
\define\End{\text{\rm End}}
\define\Aut{\text{\rm Aut}}
\define\Ind{\text{\rm Ind}}

\define\tr{\text{\rm tr}}

\define\a{\alpha}
\redefine\b{\beta}
\redefine\c{\chi}
\define\g{\gamma}

\define\e{\epsilon}

\define\io{\iota}
\redefine\o{\omega}
\define\p{\pi}
\define\ph{\phi}
\define\ps{\psi}
\define\r{\rho}
\define\s{\sigma}
\redefine\t{\tau}

\redefine\l{\lambda}
\define\z{\zeta}
\define\x{\xi}

\redefine\G{\Gamma}

\define\Om{\Omega}

\redefine\L{\Lambda}

\define\boc{\bold c}

\define\kk{\bold k}

\redefine\ss{\bold s}
\redefine\tt{\bold t}

\redefine\AA{\bold A}

\define\CC{\bold C}

\define\EE{\bold E}
\define\FF{\bold F}
\define\GG{\bold G}

\define\II{\bold I}
\define\JJ{\bold J}
\define\KK{\bold K}

\define\NN{\bold N}

\define\QQ{\bold Q}

\define\SS{\bold S}

\define\VV{\bold V}
\define\WW{\bold W}
\define\ZZ{\bold Z}
\define\XX{\bold X}

\define\ca{\Cal A}
\define\cb{\Cal B}
\define\cc{\Cal C}
\define\cd{\Cal D}
\define\ce{\Cal E}
\define\cf{\Cal F}
\define\cg{\Cal G}
\define\ch{\Cal H}

\define\cl{\Cal L}

\define\co{\Cal O}

\define\cs{\Cal S}
\define\ct{\Cal T}
\define\cu{\Cal U}
\define\cv{\Cal V}
\define\cw{\Cal W}
\define\cz{\Cal Z}
\define\cx{\Cal X}

\define\fb{\frak b}
\define\fc{\frak c}

\define\fg{\frak g}

\define\ft{\frak t}

\define\fz{\frak z}

\define\fA{\frak A}

\define\fF{\frak F}

\define\fS{\frak S}

\define\fV{\frak V}
\define\fZ{\frak Z}

\define\tg{\ti g}

\define\tj{\ti j}

\define\tq{\ti q}

\define\tu{\ti u}

\define\tM{\ti M}

\define\tW{\ti W}
\define\tX{\ti X}
\define\tY{\ti Y}

\define\sha{\sharp}

\define\tce{\ti\ce}
\define\bul{\bullet}

\define\che{\check}

\define\KL{KL}
\define\Kml{KmL1}
\define\KmL{KmL2}
\define\FOUR{L1}
\define\SING{L2}
\define\ICC{L3}
\define\CUS{L4}
\define\CUSII{L5}
\define\CLA{L6}
\define\CLAII{L7}
\define\AFF{L8}
\define\CDG{L9}
\define\CLE{L10}
\define\YU{Y}
\define\ZH{Zh}
\define\uchJ{\un{\chJ}}
\define\chJ{\che J} 
\define\chS{\che S}
\define\hcs{\hat{\cs}}
\define\bcs{\bar{\cs}}
\define\dcs{\dot{\cs}}
\define\huE{\hat{\uE}}
\define\ufZ{\un{\fZ}}
\define\hR{\hat R}
\define\uc{\un{c}}
\define\tGG{\ti{\GG}}
\define\chJJ{\che{\JJ}}
\head Introduction\endhead
\subhead 0.1\endsubhead
Let $G$ be a simple adjoint algebraic group defined and split over the finite field $\FF_q$.
Let $K_0=\FF_q((\e))$, $K=\bar{\FF_q}((\e))$. 
We are interested in the characters of the standard representations (in the
sense of Langlands) of $G(K_0)$ corresponding to the (irreducible) unipotent representations (\cite{\CLA}) 
of $G(K_0)$, restricted to the set $G(K_0)_{rsc}=G(K)_{rsc}\cap G(K_0)$
where $G(K)_{rsc}$ is the intersection of the set $G(K)_{rs}$ of regular semisimple elements in $G(K)$ with
the set $G(K)_c$ of compact elements in $G(K)$ (that is, elements which normalize some Iwahori subgroup of 
$G(K)$); we call these restrictions the {\it unipotent characters} of $G(K_0)$. 
We hope that the unipotent characters (or some small linear combination of them) have a 
geometric meaning in the same way as the characters of (irreducible) unipotent representations of 
$G(\FF_q)$ can be expressed in terms of character sheaves on $G$. Thus we are seeking some geometric objects
on $G(K)_c$ on which the Frobenius map acts and from which the unipotent characters can be recovered. 

In this paper we define a collection of class functions on $G(K_0)_{rsc}$ which we call {\it unipotent 
almost characters}. (These class functions can conceivably take the value $\iy$ at some points but we 
conjecture that the set of such points is empty; in the rest of this introduction we assume that this 
conjecture holds.) We expect that the unipotent almost characters are in the same relation with
the unipotent characters as the objects with the same names associated to $G(\FF_q)$. In particular we 
expect that the unipotent characters generate same subspace of the
vector space of class functions on $\cv_0$ as the unipotent
almost characters. (A refinement of this is stated 
as a conjecture in 3.11(a).) Moreover we expect that each unipotent almost
character can be expressed as a linear
combination of a small number of unipotent characters, just like for $G(\FF_q)$.

Our definition of unipotent almost characters is similar to one of the two definitions of the analogous
functions for $G(\FF_q)$ (which was in terms of character sheaves on $G(\bar{\FF}_q)$). They are associated 
to some new geometric objects on $G(K)_c$ which can be thought of as character sheaves on $G(K)_c$ (or 
rather, cohomology sheaves of character sheaves) and are defined even when $\kk$ is replaced by any
algebraically closed field such as $\CC$ (in which case $G(K)$ becomes $G(\CC((\epsilon))$.)

The definition of these new geometric objects combines three ingredients:

(i) A generalization of the construction \cite{\AFF} of an affine Weyl group action on the homology of the 
variety of Iwahori subgroups (see \cite{\KL}) of $G(K)$ containing a given element of $G(K)_c$.

(ii) A construction of co-standard representations of an affine Weyl group in the framework of the 
generalized Springer correspondence \cite{\ICC}.

(iii) A matching of the affine Weyl groups appearing in (i) and (ii).

Now (i) (which is discussed in \S3) is based on some preliminaries on unipotent character sheaves on 
disconnected groups given in \S1. It uses geometry (such as perverse sheaves) arising from $G(K)$. On the 
other hand, (ii) (which is discussed in \S2) is a variant of the geometric construction of representations 
of graded affine Hecke algebras given in \cite{\CUS}.

It should be pointed out that the notion of co-standard module which appears in (ii) is not intrinsic to the
affine Weyl group, but it is associated to an affine Weyl group viewed as the limit as $q\to1$ of an affine 
Hecke algebras with possibly unequal parameters $q^{m_i}$. While the irreducible representations of the 
affine Weyl group in (ii) have an elementary definition (in terms of representations of finite Weyl groups),
our definition of the co-standard representations is in terms of geometry (such as perverse sheaves) arising
from the group of type dual to that of $G$.

Since the affine Weyl groups in (i) and (ii) appear in totally different worlds (one from $G$, the other 
from the dual group of $G$) the fact that they match is a miracle (which has already been exploited in 
\cite{\CLA}).

After these new geometric objects are defined, the unipotent almost characters are defined in terms of them 
by taking traces of the Frobenius map. 

In \S4 we give some supporting evidence, based mostly on \cite{\KmL}, for the conjectures in this paper. 

We expect that a similar picture exists with 
unipotent representations replaced by representations of depth zero. In \S5 we discuss a possible 
generalization to $p$-adic groups in unequal characteristic.

\subhead 0.2. Notation\endsubhead
If $\G$ is a group then $\cz_\G$ is the centre of $\G$. If $H$ is a subgroup of $\G$, then $N_\G H$ is the 
normalizer of $H$ in $\G$. If $g\in\G$ then $Z_\G(g)$ is the centralizer of $g$ in $\G$. We denote by 
$\Irr\G$ a set of representatives for the isomorphism classes of irreducible finite dimensional 
representations of $\G$ (over $\CC$). We fix an algebraically closed field $\kk$ and a prime number $l$ 
invertible in $\kk$. Let $\bbq$ be an algebraic closure of the field of $l$-adic numbers. For an algebraic 
variety $X$ over $\kk$ let $\cd(X)$ be the bounded derived category of constructible $\bbq$ sheaves on $X$. 
By "local system" on $X$ we usually mean a $\bbq$-local system. (An exception is in \S2 where varieties and 
local systems are over $\CC$.) If $\ce$ is a local system on $X$ and $i\in\NN$ we set
$H_{-i}(X,\ce)=(H^i_c(X,\ce^\da))^\da$; generally, $()^\da$ denotes the dual of a vector space or of a local 
system. If $x\in X$, $\ce_x$ denotes the stalk of $\ce$ at $x$. If $H$ is a linear algebraic group over 
$\kk$, let $H^0$ be the identity component of $H$. If $H$ acts on $X$ let $\cd_H(X)$ be the corresponding 
equivariant derived category.

For a finite set $S$ let $|S|$ be the cardinal of $S$.

In \S3 we assume that $\bbq,\CC$ are identified as fields.

\head 1. Preliminaries on character sheaves on disconnected groups\endhead
\subhead 1.1\endsubhead
Let $G$ be an affine algebraic group over $\kk$ such that $G^0$ is reductive. For any subgroup $H$ of $G$ we
write $NH$ instead of $N_GH$. The set of subgroups of $G^0$ containing a fixed Borel subgroup of $G^0$ is of
the form $\{B_J;J\sub I\}$ ($I$ is a finite indexing set) where for $J\sub I,J'\sub I$ we have 
$B_J\sub B_{J'}$ if and only if $J\sub J'$. In particular $B_\em$ is a 
Borel subgroup. Let $W$ be a (finite)
indexing set for the set of $(B_\em,B_\em)$ double cosets in $G$. For $w\in W$ let $O_w$ be the 
corresponding double coset. Let $W'=\{w\in W;O_w\in G^0\}$ and let $\Xi=\{w\in W;O_w\sub NB_\em\}$. If 
$i\in I$ then $B_{\{i\}}-B_\em=O_w$ for a well defined $w\in W'$; we set $w=s_i$. There is a unique group 
structure on $W$ such that the following holds: if $w,w'\in W$ are such that $O_wO_{w'}$ is of the form 
$O_{w''}$ for some $w''\in W$ then $ww'=w''$; if $i\in I$ then $s_i^2=1$. Now $W'$ is the subgroup of $W$ 
generated by $\{s_i;i\in I\}$; it is a Coxeter group on these generators. Also $\Xi$ is a subgroup of $W$ 
such that $W=\Xi W'=W'\Xi$, $\Xi\cap W'=\{1\}$. If $\x\in\Xi$ and $i\in I$ then $\x s_i\x\i=s_{\x(i)}$ for a
unique $\x(i)\in I$; moreover $\x:i\m\x(i)$ is an action of $\Xi$ on $I$. For $J\sub I$, let $W_J$ be the 
subgroup of $W'$ generated by $\{s_i;i\in J\}$ and let $\Xi_J=\Xi\cap N_WW_J$; let $w_0^J$ be the longest 
element of the finite Coxeter group $W_J$. For $J\sub I$ let $\chJ=I-J$ and let ${}^{\chJ}W$ (resp. 
${}^{\chJ}W'$) be the set of all $w\in N_WW_J$ (resp. $w\in N_{W'}W_J$) such that $w$ has minimal length in 
its $W_J$-coset; this is a subgroup of $N_WW_J$ (resp. $N_{W'}W_J$).

For $\x\in\Xi$ we set ${}^\x G=\cup_{w\in W'\x}O_w$. Note that $G=\sqc_{\x\in\Xi}{}^\x G$ is the 
decomposition of $G$ into connected components. 

\subhead 1.2\endsubhead
For $J\sub I$, $B_J$ is a parabolic subgroup of $G^0$ and $NB_J=\sqc_{\x\in\Xi_J}{}^\x NB_J$ where for 
$\x\in\Xi_J$, ${}^\x NB_J:=NB_J\cap{}^\x G$ is a single $B_J$-coset. Let $U_J$ be the unipotent radical of 
$B_J$ (a normal subgroup of $NB_J$). Let $\ov{NB_J}=NB_J/U_J$, $\ov{B_J}=B_J/U_J$, and let 
$p_J:NB_J@>>>\ov{NB_J}$ be the obvious homomorphism. For any $\x\in\Xi_J$ let 
${}^\x\ov{NB_J}=p_J({}^\x NB_J)$. Now $\ov{NB_J}$ is a connected reductive group over $\kk$ whose connected
components are ${}^\x \ov{NB_J}$ for various $\x\in\Xi_J$; the identity component is $\ov{B_J}$. For 
$\x\in\Xi_J$ let ${}_\x\cz_{\ov{B_J}}=\cz_{\ov{B_J}}\cap Z_{\ov{NB_J}}(g)$ where $g$ is any element of
${}^\x\ov{NB_J}$. 

\subhead 1.3\endsubhead
We consider a triple $(J,\o,A)$ where $J\sub I$, $\o\in\Xi_J$ and $A$ is a unipotent cuspidal character 
sheaf (see \cite{\CDG, X,44.4}) on the connected component ${}^\o\ov{NB_J}$ of $\ov{NB_J}$. By 
\cite{\CDG, VI,30.2} there is a well defined subvariety $\fS\sub{}^\o\ov{NB_J}$ and a local system $\cs$ on 
$\fS$ such that for some integer $n\ge1$ invertible in $\kk$, the following hold:

$\fS$ is a single orbit of ${}_\o\cz_{\ov{B_J}}^0\T\ov{B_J}$ acting on ${}^\o\ov{NB_J}$ by 
$(z,x):g\m xz^ngx\i$;

$\cs$ is ireducible and equivariant for this action; $A|_{\fS}=\cs[\dim\fS], A|_{{}^\o\ov{NB_J}-\fS}=0$.
\nl
(We have used that $A$ is clean (see \cite{\CDG, X}, \cite{\CLE} and its references).
Let $cl(\fS)$ be the closure of $\fS$ in ${}^\o\ov{NB_J}$. 

Let ${}^{\uchJ}W'$ be the fixed point set of the automorphism of ${}^{\chJ}W'$ induced by 
$\Ad(\o):W'@>>>W'$. If $L$ is a Levi subgroup of $B_J$ then ${}^{\uchJ}W'$ can be identified with 
$\G:=\{y\in N_{G^0}L;y(NL\cap NB_J)y\i=NL\cap NB_J\}/L$. We can view $\fS$ as a subset of $NL\cap NB_J\}/L$. 
From the classification of unipotent cuspidal character sheaves \cite{\CDG, X} we see that if 
$y\in N_{G^0}L$ represents an element of $\G$ then $\Ad(y)$ preserves $\fS$ and $\cs$. Let 
$$X=\{(g,xB_J)\in{}^\o G\T G^0/B_J;x\i gx\in {}^\o NB_J,p_J(x\i gx)\in\fS\}.$$
Let $\bY$ be the set of all $g\in{}^\o G$ such that for some $x\in G^0$ we have
$x\i gx\in{}^\o NB_J,p_J(x\i gx)\in cl(\fS)$. Let $\ps:X@>>>\bY$ be the first projection. We define a local 
system $\tcs$ on $X$ by $\tcs_{(g,xB_J)}=\cs_{p_J(x\i gx)}$; this is well defined by the 
$\ov{B_J}$-equivariance of $\cs$. Using \cite{\CDG, I,5.7}, \cite{\CDG, II,7.10(a)} we see 
that the vector space $\EE:=\End(\ps_!\tcs)$ has a canonical direct sum decomposition into lines $\EE_w$
indexed by the elements $w\in\G$ hence by the elements $w\in{}^{\uchJ}W'$. Moreover in the algebra 
structure of $\EE$ we have $\EE_w\EE_{w'}=\EE_{ww'}$ for $w,w'\in {}^{\uchJ}W'$. In particular all nonzero 
elements in $\EE_w$ are units in $\EE$.

\subhead 1.4\endsubhead
Let $\boc$ be the two-sided cell of $W_J$ associated in \cite{\CDG, X,44.18} to the unipotent character 
sheaf $A$ on ${}^\o\ov{NB_J}$. Let $\boc'$ be the two-sided cell of $W'$ that contains $\boc$. Consider the 
bijection $j$ in \cite{\CDG, X,44.21(h)} from the set of unipotent character sheaves on ${}^\o\ov{NB_J}$ 
with associated two-sided cell $\boc$ to the set of unipotent character sheaves of ${}^\o G$ with associated 
two-sided cell $\boc'$. (The assumptions of {\it loc.cit.} are satisfied by the classification of unipotent 
character sheaves.) From the definition we see that $j(A)$ appears with multiplicity $1$ in the semisimple 
perverse sheaf $\ps_!\tcs[\dim X]$ on $\bY$. Now $\EE$ acts naturally on the one-dimensional vector space
$V:=\Hom(j(A),\ps_!\tcs[\dim X])$. For $w\in{}^{\uchJ}W'$, since the nonzero elements of $\EE_w$ are units 
in $\EE$, we see that precisely one element $b_w\in\EE_w$ acts on $V$ as identity. Clearly we have 
$b_wb_{w'}=b_{ww'}$ for $w,w'\in{}^{\uchJ}W'$. Hence the basis $(b_w)$ of $\EE$ identifies the algebra $\EE$
with the group algebra $\bbq[{}^{\uchJ}W']$. Thus we have
$$\End_{\cd(\bY)}(\ps_!\tcs)=\bbq[{}^{\uchJ}W'].\tag a$$

\head 2. Co-standard representations of certain (extended) affine Weyl groups\endhead
\subhead 2.1\endsubhead
Let $\cg$ be a connected, almost simple, simply connected algebraic group over $\CC$. We fix a maximal torus
$T$ of $\cg$ with Lie algebra $\ft$. Let $R\sub\ft^\da$ be the set of roots. Let $\a_i\in R$ ($i\in[0,n]$) be
such that $(\a_i)_{i\in[1,n]}$ form a set of simple roots for $R$ and $\a_0$ is the negative of the 
corresponding highest root. For $i\in[0,n]$ let $h_i\in\ft$ be the coroot corresponding to $\a_i$. There are
unique integers $n_i>0$ $(i\in[0,n])$ such that $\sum_{i\in[0,n]}n_i\a_i=0$, $n_0=1$. Let $V$ be a 
$\CC$-vector space with basis $\{b_i;i\in[0,n]\}$ and let $\{b'_i;i\in[0,n]\}$ be the dual basis of $V^\da$.
The canonical pairing $V\T V^\da@>>>\CC$ is denoted by $x,x'\m x(x')$. We imbed $\ft$ into $V^\da$ by 
$y\m\sum_{i\in[0,n]}\a_i(y)b'_i$; we identify $\ft$ with its image, the subspace 
$\{\sum_{i\in[0,n]}c_ib'_i;c_i\in\CC,\sum_in_ic_i=0\}$. In particular we regard $h_i$ as a vector in 
$V^\da$. For any $y\in\ft$, $i\in[0,n]$ we have $b_i(y)=\a_i(y)$. Let
$$\ft^1=\{\sum_{i\in[0,n]}c_ib'_i;c_i\in\CC,\sum_in_ic_i=1\}\sub V^\da.$$
For $i\in[0,n]$ we define $s_i:V@>>>V$ by $s_i(x)=x-x(h_i)b_i$ and its contragredient $s_i:V^\da@>>>V^\da$ 
by $s_i(x')=x'-b_i(x')h_i$.
Let $W$ be the subgroup of $GL(V)$ or $GL(V^\da)$ generated by $\{s_i;i\in[0,n]\}$ (an affine Weyl group).
Note that $\ft,\ft^1$ are $W$-stable subsets of $V^\da$.
For any $J\sub[0,n]$ we set $\chJ=[0,n]-J$. 

For any $S\sub[0,n],S\ne\em$ let 
$$C_S=\{x'\in\ft^1;x'=\sum_{i\in S}c_ib'_i\text{ with }c_i\in\CC,c_i>0\qua\frl i\in S\}.$$
Here, for a complex number $c=a+b\sqrt{-1}$ with $a,b$ real, we write $c>0$ whenever $a>0$ or $a=0,b>0$. 
The sets $C_S$ are disjoint.
For any $J\sneq[0,n]$ let $W_J$ be the subgroup of $W$ generated by $\{s_i;i\in J\}$ (a finite Coxeter group
whose longest element is denoted by $w_0^J$). Let ${}^{\chJ}W$ be the set of all $w\in N_WW_J$ such that $w$
has minimal length in its $W_J$-coset; this is a subgroup of $N_WW_J$.

For any $J\sneq[0,n]$ let $R_J$ (resp. $R_J^+$) be the set of $\a\in R$ which are $\ZZ$-linear (resp. 
$\NN$-linear) combinations of $\{a_i;i\in J\}$ and let $\cg_J$ be the subgroup of $\cg$ generated by $T$ and
by the root subgroups of $\cg$ corresponding to the roots in $R_J$; note that $\cg_J$ is a connected
reductive subgroup of $\cg$ with maximal torus $T$ and root system $R_J$. Let $\fz_J$ be the Lie algebra of 
the centre of $\cg_J$. Let $V_J$ be the subspace of $V$ spanned by $\{b_i;i\in J\}$.
Let $V^\da_{\chJ}$ be the subspace of $V^\da$ spanned by $\{b'_i;i\in\chJ\}$. We have 
$\fz_J=\ft\cap V^\da_{\chJ}$. Let $\fz_J^1=\ft^1\cap V^\da_{\chJ}$, an affine hypersurface in $\fz_J$.

\subhead 2.2\endsubhead
Let $\fZ_\cg$ be the set of all pairs $(\fc,\fF)$ where $\fc$ is a conjugacy class in $\cg$ and $\fF$ is an
irreducible $\cg$-equivariant local system on $\fc$ (up to isomorphism). We have a partition 
$$\fZ_\cg=\sqc_{\c\in\Hom(\cz_\cg,\CC^*)}\fZ_{\cg,\c}$$ 
where $\fZ_{\cg,\c}$ is the set of all $(\fc,\fF)\in\fZ_\cg$ such that the character by which $\cz_\cg$ acts
on any stalk of $\fF$ in the $\cg$-equivariant structure of $\fF$ is $\c$.

\subhead 2.3\endsubhead
We fix a triple $(J,\cc,\ce)$ where $J\sneq[0,n]$, $\cc$ is a unipotent conjugacy class of $\cg_J$ and $\ce$
is an irreducible cuspidal $\cg_J$-equivariant local system on $\cc$. From the classification of cuspidal 
local systems \cite{\ICC} we see that $J$ has the following property:

(a) {\it if $|\chJ|\ge2$ then for any $k\in\chJ$ we have 
$\ss_k:=w_0^{J\cup k}w_0^J=w_0^Jw_0^{J\cup k}\in{}^{\chJ}W$.}
\nl
From \cite{\CLA, 1.15} we see that the following holds.

(b) {\it If $|\chJ|\ge2$ then ${}^{\chJ}W$ is a Coxeter group (an irreducible affine Weyl group) on
 generators $\ss_k$ ($k\in\chJ$).}
\nl
On the other hand, if $|\chJ|=1$ we have ${}^{\chJ}W=\{1\}$. 

\subhead 2.4\endsubhead
We preserve the setup of 2.3. Assume that $|\chJ|\ge2$. Then each generator $\ss_k$ of ${}^{\chJ}W$ leaves 
stable the subspace $V_J$ of $V$ hence also its annihilator $V^\da_{\chJ}$ in $V^\da$. Hence ${}^{\chJ}W$ 
acts on $V^\da_{\chJ}$. This action leaves stable the subset $\fz^1_J$ of $V^\da_{\chJ}$ and also the 
subspace $\fz_J$ of $V^\da_{\chJ}$ where it acts through a finite quotient $\cw_{\chJ}\sub GL(\fz_J)$.

Let $\cl'$ be the set of all $x\in\fz_J$ such that the translation $z\m z+x$ of $\fz^1_J$ coincides with the
automorphism $x\m w(x)$ of $\fz^1_J$ for some $w\in{}^{\chJ}W$. Then $\cl'$ is a subgroup of $\fz_J$ such 
that $\CC\ot\cl'=\fz_J$. Let $\cl=\{x\in\fz_J^\da;x(\cl')\sub\ZZ\}$. Clearly $\cw_{\chJ}$ acts naturally on 
$\cl',\cl$. Hence we can form the semidirect product 
$$\tW_{\chJ}=\cl\cdot\cw_{\chJ}$$
(with $\cl$ normal). Let 
$$\ct=\CC^*\ot\cl'=\Hom(\cl,\CC^*).$$ 
For any $t\in\ct$ we denote by $\cw_{\chJ,t}$ the stabilizer of $t$ in $\cw_{\chJ}$ (for the action on 
$\ct$) and we form the semidirect product $\tW_{\chJ,t}=\cl\cdot\cw_{\chJ,t}$ (with $\cl$ normal); this is a
subgroup of finite index of $\tW_{\chJ}$.

Let $\exp_\ct:\fz_J@>>>\ct$ (that is $\CC\ot\cl'@>>>\CC^*\ot\cl$) be the homomorphism induced by the 
exponential $\CC@>>>\CC^*$. By \cite{\CLAII, 9.2(b)} we can find $k_J\in\chJ$ so that the following holds.

(a) {\it Let $p_J:\fz^1_J@>>>\ct$ be the composition of the map $x\m x-(1/n_{k_J})b'_{k_J}$, 
$\fz^1_J@>>>\fz_J$ with the map $x\m\exp_\ct(2\pi\sqrt{-1}x)$, $\fz_J@>>>\ct$. Then $x'\m p_J(x')$ is a 
bijection between $D_J:=\cup_{S\sub\chJ;S\ne\em}C_S$ and a set of representatives for the 
$\cw_{\chJ}$-orbits in $\ct$.}
\nl
(In the case where $J=\em$ we take $k_\em=0$.) By Mackey's theorem, $\Irr\tW_{\chJ}$ is in natural bijection
with the set of pairs $(t,\r)$ where $t$ runs over a set of representatives for the $\cw_{\chJ}$-orbits on 
$\Hom(\cl,\CC^*)=\ct$ and $\r\in\Irr\cw_{J,t}$. Using (a) this can be viewed as a bijection
$$\Irr\tW_{\chJ}\lra\{(d,\r);d\in D_J,\r\in\Irr\cw_{\chJ,p_J(d)}\}.\tag b$$
From \cite{\CLA, 3.8,3.9} we deduce:

(c) {\it If $d\in C_S$ (with $S\sub\chJ,S\ne\em$) then the subgroup $({}^{\chJ}W)_{\chS-J}$ of ${}^{\chJ}W$ 
generated by $\{\ss_k;k\in\chS-J\}$ maps isomorphically to $\cw_{\chJ,p'_J(d)}$ under the canonical 
homomorphism ${}^{\chJ}W@>>>\cw_{\chJ}$.}
\nl
Now the natural action of $({}^{\chJ}W)_{\chS-J}$ on $\cl$ (restriction of the ${}^{\chJ}W$-action) is the 
same as the action of $\cw_{\chJ,p_J(d)}$ in (c) on $\cl$. Using this action we can form the semidirect 
product $\tW_{\chJ,\chS-J}=\cl\cdot({}^{\chJ}W)_{\chS-J}$ (with $\cl$ normal); this is a subgroup of finite 
index of $\tW_{\chJ}$ (using the inclusion $({}^{\chJ}W)_{\chS-J}\sub\cw_{\chJ}$ coming from (c)).

Hence (b) can be viewed as a bijection
$$\Irr\tW_{\chJ}\lra\sqc_{S\sub\chJ;S\ne\em}\{(d,\r);d\in C_S,\r\in\Irr({}^{\chJ}W)_{\chS-J}\}\tag d$$ 
which associates to $(d,\r)$ (with $d\in C_S$) the irreducible representation 
$$\Ind_{\tW_{\chJ,\chS-J}}^{\tW_{\chJ}}([p_J(d)]\ot\ti\r)\tag e$$
where $\ti\r$ is the irreducible representation of $\tW_{\chJ,\chS-J}$ which equals $\r$ on 
$({}^{\chJ}W)_{\chS-J}$ and on which $\cl$ acts as identity; $[p_J(d)]$ is the one dimensional 
representation of $\tW_{\chJ,\chS-J}$ which equals $p_J(d)$ on $\cl$ (recall that 
$p_J(d)\in\Hom(\cl,\CC^*)$) and on which $({}^{\chJ}W)_{\chS-J}$ acts trivially.

\subhead 2.5\endsubhead
We preserve the setup of 2.3. Assume that $|\chJ|=1$. In this case we set $\tW_{\chJ}=\{1\}$ and 
$({}^{\chJ}W)_\em=\{1\}$. Then a bijection as in 2.4(d) continues to hold (both sides have exactly one 
element; the only $S$ in the union is $S=\chJ$ and the corresponding $C_S$ has exactly one element.)

\subhead 2.6\endsubhead
Putting together the bijections 2.4(e) (see also 2.5) for various $(J,\cc,\ce)$ as in 2.3 (with $\ce$ defined
up to isomorphism) we obtain identifications
$$\align&\sqc_{(J,\cc,\ce)}\Irr\tW_{\chJ}=
\sqc_{(J,\cc,\ce)}\sqc_{S\sub\che J;S\ne\em}\{(d,\r);d\in C_S,\r\in\Irr({}^{\chJ}W)_{\chS-J}\}\\&
=\sqc_{S\sub[0,n];S\ne\em}\sqc_{d\in C_S}\sqc_{(J,\cc,\ce);J\sub\chS}\Irr({}^{\chJ}W)_{\chS-J}.\tag a
\endalign$$
We can identify $({}^{\chJ}W)_{\chS-J}$ in an obvious way with the group ${}^{\chS-J}(W_{\chS})$ of all 
$w\in N_{W_{\chS}}(W_{\chS-J})$ which have minimal length in their $W_{\chS-J}$-coset.

By the generalized Springer correspondence \cite{\ICC, \S6,\S9} applied to the connected reductive group 
$\cg_{\chS}$ we have for any $S\sub[0,n],S\ne\em$ an identification of
$$\sqc_{(J,\cc,\ce);J\sub\chS}\Irr({}^{\chJ}W)_{\chS-J}=
\sqc_{(J,\cc,\ce);J\sub\chS}\Irr{}^{\chS-J}(W_{\chS})$$
with the set $\ca_{\chS}$ consisting of all pairs $(c,\cf)$ where $c$ is a unipotent class of $\cg_{\chS}$ 
and $\cf$ is an irreducible $\cg_{\chS}$-equivariant local system on $c$ (up to isomorphism). We shall 
denote by $\ca_{\chS,J,\cc,\ce}$ the subset of $\ca_{\chS}$ corresponding to a fixed $(J,\cc,\ce)$ 
under the previous identification.

Introducing this into (a) we obtain an identification
$$\sqc_{(J,\cc,\ce)}\Irr\tW_{\chJ}=\sqc_{S\sub[0,n];S\ne\em}\sqc_{d\in C_S}\ca_{\chS}.\tag b$$
Applying 2.4(a) with $J=\em$ (so that $\ct=T$) we see that $p_\em:\ft^1@>>>T$ defines a bijection between 
$\cup_{S\sub[0,n];S\ne\em}C_S$ and a set of representatives for the Weyl group orbits in $T$. Note also that
if $d\in C_S$ with $S\sub[0,n],S\ne\em$ then $Z_\cg(p_\em(d))=\cg_{\chS}$. Hence associating to $S,d$ as in
(b) and to $(c,\cf)\in\ca_{\chS}$ the pair $(\fc,\cf')$ where $\fc$ is the conjugacy class in $\cg$ 
containing $p_\em(d)c$ and $\cf'$ is the irreducible $\cg$-equivariant local system on $\fc$ whose 
restriction to $p_\em(d)c\cong c$ is $\cf$, gives a bijection between the right hand side of (b) and 
$\fZ_\cg$ (see 2.2). Now (b) becomes a bijection 
$$\sqc_{(J,\cc,\ce)}\Irr\tW_{\chJ}\lra\fZ_\cg.\tag c$$

\subhead 2.7\endsubhead
Let $(J,\cc,\ce)$ be as in 2.3. Let $E\in\Irr\tW_{\chJ}$. We write $E=E_{\fc,\fF}$ where 
$(\fc,\fF)\in\fZ_\cg$ corresponds to $E$ under 2.6(c). We also write $E=E^{S,d,c,\cf}$ where $S\sub\chJ$,
$S\ne\em$, $d\in C_S$ and $(c,\cf)\in\ca_{\chS,J,\cc,\ce}$ (see 2.4) are such that $E$ corresponds to 
$S,d,c,\cf$ under 2.6(b).

We write $E$ in the form $E=\Ind_{\tW_{\chJ,\chS-J}}^{\tW_{\chJ}}([p_J(d)]\ot\ti\r)$ as in 2.4(a). Recall 
that $\r$ is an irreducible representation of $({}^{\chJ}W)_{\chS-J}={}^{\chS-J}(W_{\chS})$ and $\ti\r$ is 
the irreducible representation of $\tW_{\chJ,\chS-J}$ which equals $\r$ on $({}^{\chJ}W)_{\chS-J}$ and on
which $\cl$ acts as identity. 

Our next objective is to define a $\tW_{\chJ,\chS-J}$-module $\hat{\ti\r}$ (of finite dimension) which 
contains $\ti\r$ as a $\tW_{\chJ,\chS-J}$-submodule.

Let $P$ be the parabolic subgroup of $\cg_{\chS}$ which has $\cg_J$ as a Levi subgroup and contains the 
root subgroups corresponding to roots in $R_{\chS}^+$; let $U_P$ be the unipotent radical of $P$. Let 
$g\in c$, let $\fA=Z_{\cg_{\chS}}(g)/Z_{\cg_{\chS}}^0(g)$ and let $\t$ be the irreducible representation of 
$\fA$ corresponding to the local system $\cf$ on $c$. Let 
$$X=\{xP\in\cg_{\che S}/P;x\i gx\in\cc U_P\}.$$
By \cite{\ICC, 1.2(b)}, we have $\dim X\le e/2$ where 
$$e=\dim\cg_{\chS}-\dim\cg_J-\dim c+\dim\cc.$$
Let $\tce$ be the local system on $X$ such that $\tce_{xP}$ is the stalk of $\ce$ at the $\cc$-component of 
$x\i gx\in\cc U_P$. By \cite{\ICC, \S9} for any $i\ge0$ there is a natural degree preserving action of 
${}^{\chS-J}(W_{\chS})$ on $H^\bul(X,\tce)=\op_{i\ge0}H^i_c(X,\tce)$ which commutes with the obvious action
of $\fA$. (By \cite{\CUS, 8.6} we have $H^i_c(X,\tce)=0$ unless $i$ is even.) Moreover, the induced 
$({}^{\chJ}W)_{\chS-J}$-action on $\Hom_{\fA}(\t,H^e_c(X,\tce))$ is the same as $\r$. In particular, we have
$\dim X=e/2$. Let $\SS=\op_{r\ge0}\SS_{2r}$ be the symmetric algebra of $\fz_J^\da$ with its grading in 
which the elements of $\fz_J^\da$ have degree $2$. We can identify $\SS$ with the total 
$\cg_{\chS}$-equivariant cohomology algebra $H^\bul_{\cg_{\chS}}(X,\bbq)$ of $X$ (for the natural action of 
$\cg_{\chS}$) by arguing as in \cite{\CUS, 4.2}: after choosing $g_0\in\cc$ and setting $M=Z_{\cg_J}(g_0)$, 
we have
$$H^\bul_{\cg_{\chS}}(X,\bbq)=H^\bul_P(\cc U_P)=H^\bul_P(\cc)=H^\bul_{\cg_J}(\cc)
=H^\bul_M(\text{point})=H^\bul_{\cz_{\cg_J}^0}(\text{point})=\SS.$$
(We use that $M^0$ is $\cz_{\cg_J}^0$ times a unipotent group; $P,\cg_J$ act by conjugation.)
Now we have an obvious algebra homomorphism $H^\bul_{\cg_{\chS}}(X,\bbq)@>>>H^\bul_{\{1\}}(X,\bbq)$ (see 
\cite{\CUS, 1.4(g)}) or equivalently $\SS@>>>H^\bul(X,\bbq)$. Note also that 
$H^\bul_c(X,\tce)$ is naturally an $H^\bul(X,\bbq)$-module (via the cup product) hence via the previous 
algebra homomorphism, $H^\bul_c(X,\tce)$ is naturally an $\SS$-module. This module structure combines with 
the $({}^{\chJ}W)_{\chS-J}$-module structure on $H^\bul_c(X,\tce)$ to give a module structure on 
$H^\bul_c(X,\tce)$ over $\SS\ot\bbq[({}^{\chJ}W)_{\chS-J}]$ which is regarded as a $\bbq$-algebra in which 
$(x\ot w)(x'\ot w')=(xw(x'))\ot(ww')$ for $x,x'\in\SS$ and $w,w'\in({}^{\chJ}W)_{\chS-J}$ (we use the action
of $({}^{\chJ}W)_{\chS-J}$ on $\SS$ coming from its natural action on $\fz_J^\da$, restriction of the 
${}^{\chJ}W$-action); this follows by specialization from \cite{\CUS, 8.13}. In the $\SS$-module structure 
on $H^\bul_c(X,\tce)$ we have 
$$\SS_{2r}H^i_c(X,\tce^\da)\sub H^{i+2r}_c(X,\tce^\da)\tag a$$
for any $i\ge0$. Let $\hat\SS$ be the completion of $\SS$ at the maximal ideal $\op_{r>0}\SS_{2r}$. Since 
$H^i_c(X,\tce)=0$ unless $i\le e$ (we use that $\dim X\le e/2$) it follows that 
$\SS_{e+2}H^\bul_c(X,\tce^\da)=0$; hence the $\SS$-module structure on $H^\bul_c(X,\tce^\da)$ extends to an 
$\hat\SS$-module structure which combines with the $({}^{\chJ}W)_{\chS-J}$-module structure on 
$H^\bul_c(X,\tce^\da)$ to give a module structure on $H^\bul_c(X,\tce^\da)$ over 
$\hat\SS\ot\bbq[({}^{\chJ}W)_{\chS-J}]$ which is regarded as a $\bbq$-algebra in the same way as 
$\SS\ot\bbq[({}^{\chJ}W)_{\chS-J}]$. Now $\tW_{\chJ,\chS-J}=\cl\cdot({}^{\chJ}W)_{\chS-J}$ can be regarded 
as a subgroup of the group of units of the algebra $\hat\SS\ot\bbq[({}^{\chJ}W)_{\chS-J}]$ (for 
$x\in\cl,w\in({}^{\chJ})W_{\chS-J}$ we view $x$ as an element of $\fz_J^\da\sub\SS$ and we associate to 
$xw\in\cl\cdot({}^{\chJ}W)_{\chS-J}$ the element $(1+x+x^2/2!+x^3/3!+\do)\ot w$ of 
$\hat\SS\ot\bbq[({}^{\chJ}W)_{\chS-J}]$). Hence the $\hat\SS\ot\bbq[({}^{\chJ}W)_{\chS-J}]$-module 
structure on $H^\bul_c(X,\tce^\da)$ restricts to a $\tW_{\chJ,\chS-J}$-module structure on 
$H^\bul_c(X,\tce^\da)$.

We set $\hat{\ti\r}=\Hom_{\fA}(\t,H^\bul_c(X,\tce^\da))$. Note that the $\fA$ action on 
$H^\bul_c(X,\tce^\da)$ commutes with the $\tW_{\chJ,\chS-J}$-module structure hence $\hat{\ti\r}$ inherits a
$\tW_{\chJ,\chS-J}$-module structure from that on $H^\bul_c(X,\tce^\da)$. From (a) we see that for any 
$a\in\ZZ$, $\op_{i;i\ge a}H^i_c(X,\tce^\da)$ is a $\tW_{\chJ,\chS-J}$-submodule of $H^\bul_c(X,\tce^\da)$. 
In particular, $H^e_c(X,\tce^\da)$ is a $\tW_{\chJ,\chS-J}$-submodule of $H^\bul_c(X,\tce^\da)$. (We use 
that $\dim X\le e/2$.) It follows that for any $a\in\ZZ$, 
$$\hat{\ti\r}_{\ge a}:=\Hom_{\fA}(\t,\op_{i;i\ge a}H^i_c(X,\tce^\da))$$ 
is a $\tW_{\chJ,\chS-J}$-submodule of $\hat{\ti\r}$. Note that $\hat{\ti\r}_{\ge e}=\ti\r$ as a
$\tW_{\chJ,\chS-J}$-module. Moreover for any $a<e$, the $\tW_{\chJ,\chS-J}$-module 
$\hat{\ti\r}_{\ge a}/\hat{\ti\r}_{\ge a-1}$ is a direct sum of $\tW_{\chJ,\chS-J}$-modules on which $\cl$ 
acts trivially and $({}^{\chJ}W)_{\chS-J}$ acts by an irreducible representation which corresponds under the 
generalized Springer correspondence to a pair $(c',\cf')\in\ca_{\chS,J,\cc,\ce}$ with $c\sub cl(c')-c'$ 
where $cl(c')$ is the closure of $c'$. We now set 
$$\hE=\Ind_{\tW_{\chJ,\chS-J}}^{\tW_{\chJ}}([p'_J(d)]\ot\hat{\ti\r}).$$
For any $a\in\ZZ$ we set 
$$\hE_{\ge a}=\Ind_{\tW_{\chJ,\chS-J}}^{\tW_{\chJ}}([p'_J(d)]\ot\hat{\ti\r}_{\ge a}).$$
Note that $\hE$ is a $\tW_{\chJ}$-module of finite dimension and that
$\do\sub\hE_{\ge a}\sub\hE_{\ge a-1}\sub\do$ are $\tW_{\chJ}$-submodules of $\hE$ such that 
$\hE_{\ge0}=\hE$, $\hE_{\ge e+1}=0$. Moreover, $\hE_{\ge e}=E$ as a $\tW_{\chJ}$-module and for any $a<e$, 
the $\tW_{\chJ}$-module $\hE_{\ge a}/\hE_{\ge a-1}$ is a direct sum of irreducible representations of 
$\tW_{\chJ}$ of the form $E^{S,d;c',\cf'}$ where $(c',\cf')\in\ca_{\chS,J,\cc,\ce}$ with $c\sub cl(c')-c'$.

We say that $\hE=\hE_{\fc,\fF}$ is the {\it co-standard} $\tW_{\chJ}$-module associated to $E=E_{\fc,\fF}$
with $(\fc,\fF)\in\fZ_\cg$.

\head 3. The main construction\endhead
\subhead 3.1 \endsubhead
Let $K=\kk((\e))$ where $\e$ is an indeterminate. Let $G$ be a simple adjoint algebraic group over $\kk$. 
Let $\GG=G(K)$. Then the parahoric subgroups of $\GG$ are well defined. For any subgroup $H$ of $\GG$ we 
write $NH$ instead of $N_\GG H$. The set of parahoric subgroups of $\GG$ containing a fixed 
Iwahori subgroup of $\GG$ is of the form $\{P_\JJ;\JJ\sneq\II\}$ ($\II$ is a finite indexing set) where for 
$\JJ\sneq\II$, $\JJ'\sneq\II$ we have $P_\JJ\sub P_{\JJ'}$ if and only if $\JJ\sub\JJ'$. In particular 
$P_\em$ is an Iwahori subgroup. 

Let $\WW$ be an indexing set for the set of $(P_\em,P_\em)$ double cosets in $\GG$. For $w\in\WW$ let $O_w$ 
be the corresponding double coset. Let $\Om=\{w\in\WW;O_w\sub NP_\em\}$. 

If $i\in\II$ then $P_{\{i\}}-P_\em=O_w$ for a well defined $w\in\WW$; we 
set $w=s_i$. We regard $\WW$ as a group as in \cite{\CLA, 1.10}; let $\WW'$ be the subgroup of $\WW$ 
generated by $\{s_i;i\in\II\}$. Note that $\WW'$ is normal in $\WW$, $\Om$ is a finite abelian subgroup of 
$\WW$ such that $\WW=\Om\WW'=\WW'\Om$, $\Om\cap\WW'=\{1\}$. Also $\WW'$ together with $\{s_i;i\in\II\}$ is a
Coxeter group (an irreducible affine Weyl group). If $\x\in\Om$ and $i\in\II$ then $\x s_i\x\i=s_{\x(i)}$ 
for a unique $\x(i)\in\II$; moreover $\x:i\m\x(i)$ is an action of $\Om$ on $\II$.

For $\JJ\sub\II$, let $\WW_\JJ$ be the subgroup of $\WW'$ generated by $\{s_i;i\in\JJ\}$ and let
$\Om_\JJ=\Om\cap N_\WW\WW_\JJ$. For $\JJ\sneq\II$ let $w_0^\JJ$ be the longest element of the finite Coxeter
group $\WW_\JJ$.

For $\JJ\sneq\II$ let $\chJJ=\II-\JJ$ and let ${}^{\chJJ}\WW$ be the set of all $w\in N_\WW\WW_\JJ$ such 
that $w$ has minimal length in its $\WW_\JJ$-coset; this is a subgroup of $N_\WW\WW_\JJ$. Let 
${}^{\chJJ}\WW'={}^{\ch\JJ}\WW\cap\WW'$. For $\JJ\sub\JJ'\sneq\II$ let 
${}^{\JJ'-\JJ}(\WW_{\JJ'})={}^{\chJJ}\WW\cap\WW_{\JJ'}$.

Let $\GG_c$ be the set of elements of $\GG$ which normalize some Iwahori subgroup.
Let $\GG_{rsc}$ be the set of regular semisimple elements in $\GG_c$.
For $\x\in\Om$ we set ${}^\x\GG=\cup_{w\in\x\WW'}O_w$, 
${}^\x\GG_c={}^\x\GG\cap\GG_c$,
${}^\x\GG_{rsc}={}^\x\GG\cap\GG_{rsc}$. Note that 
$\GG=\sqc_{\x\in\Om}{}^\x\GG$ is the decomposition of $\GG$ into left (or right) ${}^1\GG$-cosets and each 
${}^\x\GG$ is stable under $\GG$-conjugacy. 

For $\JJ\sneq\II$, $P_\JJ$ is a parahoric subgroup of $\GG$ and $NP_\JJ=\sqc_{\x\in\Om_J}{}^\x NP_\JJ$ where
for $\x\in\Om_\JJ$, ${}^\x NP_\JJ:=NP_\JJ\cap{}^\x\GG$ is a single $P_\JJ$-coset. Now $P_\JJ$ has a 
prounipotent radical $U_\JJ$ (a normal subgroup of $NP_\JJ$). Let $\ov{NP_\JJ}=NP_\JJ/U_\JJ$, 
$\ov{P}_\JJ=P_\JJ/U_\JJ$ and let $p_\JJ:NP_\JJ@>>>\ov{NP}_\JJ$ be the obvious homomorphism. For any 
$\x\in\Om_\JJ$ we set ${}^\x\ov{NP}_\JJ=p_J({}^\x NP_\JJ)$. Now $\ov{NP}_\JJ$ is naturally a reductive group
over $\kk$ whose connected components are ${}^\x\ov{NP}_\JJ$ for various $\x\in\Om_\JJ$; the identity 
component is $\ov{P}_\JJ$. For $\x\in\Om_\JJ$ let 
$${}_\x\cz_{\ov{P}_\JJ}=\cz_{\ov{P}_\JJ}\cap Z_{\ov{NP}_\JJ}(g)$$ 
where $g$ is any element of ${}^\x\ov{NP}_\JJ$.

Let $\cb$ be the set of Iwahori subgroups of $\GG$. For any $\x\in\Om$ let 
${}^\x\cx=\{(g,B)\in{}^\x\GG_c\T\cb;gBg\i=B\}$. We show that on ${}^\x\cx$ there is a notion of $l$-adic 
constructible sheaf. For simplicity we assume that $\x=1$. We have
${}^1\cx=\{(g,B)\in{}^1\GG_c\T\cb;g\in B\}$. 
For any $B\in\cb$ let $B=B_0\supset B_1\supset B_2\supset\do$ be 
the Moy-Prasad filtration of $B$; note that each $B_n$ is a normal subgroup of $B$, $B/B_n$ is an algebraic 
group over $\kk$ and $B=\underset\gets\to\lim B/B_n$. For any $n\ge1$ let 
$\cx_n=\{(g,B)\in\cx;g\in B-B_n\}$. For any $m\ge n\ge1$ let 
$\cx_{m,n}=\{(gB_m,B);B\in\cb,gB_m\in(B-B_n)/B_m\}$. (Note that $B-B_n$ is stable under left and right 
translation by $B_m$.) We define $\p_{m,n}:\cx_n@>>>\cx_{m,n}$ by $(g,B)\m(gB_m,B)$ and  
$\p'_{m,n}:\cx_{m,n}@>>>\cb$ by $(gB_m,B)\m B$. Note that the fibre of $\p'_{m,n}$ at $B\in\cb$ is the 
algebraic variety $(B-B_n)/B_m$ over $\kk$ (an open subvariety of $B/B_m$). Now $\cb$ is an inductive limit 
of projective varieties over $\kk$; taking inverse images of these projective (sub)varieties of $\cb$ under 
$\p'_{m,n}$ we obtain algebraic varieties over $\kk$ of which $\cx_{m,n}$ is the inductive limit. Hence on 
each $\cx_{m,n}$ we have a well defined notion of $l$-adic constructible sheaf. By definition, an $l$-adic 
constructible sheaf on $\cx_n$ is the inverse image of an $l$-adic constructible sheaf on $\cx_{m,n}$ under 
$\p_{m,n}$ for some $m$ such that $m\ge n$. Now an $l$-adic constructible sheaf on $\cx$ is a collection of 
$l$-adic constructible sheaves $\cf_n$ on $\cx_n$ for various $n\ge1$ such that for any $n\le n'$ the 
restriction of $\cf_{n'}$ to $\cx_n$ is $\cf_n$. (Note that $\cx_1\sub\cx_2\sub\do$ and
$\cx=\cup_{n\ge1}\cx_n$.) 

\subhead 3.2\endsubhead
Let $\fV$ be the set of all triples $(\JJ,\o,A)$ where $\JJ$ runs over a set of representatives for the
$\Om$-orbits of proper subsets of $\II$, $\o\in\Om_\JJ$ and $A$ is a unipotent cuspidal character sheaf 
(defined up to isomorphism) on the connected component ${}^\o\ov{NP}_\JJ$ of $\ov{NP}_\JJ$.

We now fix $(\JJ,\o,A)\in\fV$. By \cite{\CDG, VI,30.2} there is a well defined subvariety 
$\fS\sub{}^\o\ov{NP}_\JJ$ and a local system $\cs$ on $\fS$ such that for some integer $n\ge1$ invertible in
$\kk$, the following hold:

$\fS$ is a single orbit of $({}_\o\cz_{\ov{P}_\JJ})^0\T\ov{P}_\JJ$ acting on ${}^\o\ov{NP}_\JJ$ by 
$(z,x):g\m xz^ngx\i$;

$\cs$ is ireducible and equivariant for this action; $A|_{\fS}=\cs[\dim\fS], A|_{{}^\o\ov{NP}_\JJ-\fS}=0$.
\nl
Let $cl(\fS)$ be the closure of $\fS$ in ${}^\o\ov{NP}_\JJ$. 

Let ${}^{\un{\chJJ}}\WW=\{w\in{}^{\chJJ}\WW,\o w\o\i=w\}$,
${}^{\un{\chJJ}}\WW'=\{w\in{}^{\chJJ}\WW',\o w\o\i=w\}$. For $\JJ\sub\JJ'\sneq\II$ with $\o(\JJ')=\JJ'$,
let ${}^{\un{\JJ'-\JJ}}(\WW_{\JJ'})=\{w\in{}^{\JJ'-\JJ}(\WW_{\JJ'});\o w\o\i=w\}$.

Now $\ov{NP}_\JJ$ acts by conjugation on ${}^\o\ov{NP}_\JJ$. This action preserves $\fS$; moreover $\cs$ 
admits a $\ov{NP}_\JJ$-equivariant structure. (This can be seen from the classification of triples 
$(\JJ,\o,A)$ as above which is the same as that of the "arithmetic diagrams" in \cite{\CLA, \S7}.) We assume
that such an equivariant structure has been chosen. For any $\x\in\Om_\JJ$ and any $\s\in\fS$ the vector 
spaces $\{\cs_{z\s z\i};z\in{}^\x\ov{NP}_\JJ\}$ are canonically isomorphic to each other (by the 
$\ov{P}_\JJ$-equivariance of $\cs$) hence they can be identified with a single vector space denoted by 
$\cs^\x_\s$; the vector spaces $\cs^\x_\s$ for various $\s\in\fS$ form the stalks of a local system $\cs^\x$ 
on $\fS$. The $\ov{NP}_\JJ$-equivariant structure on $\cs$ provides an isomorphism of local systems 
$\cs@>\si>>\cs^\x$ on $\fS$.

\subhead 3.3\endsubhead
We preserve the setup of 3.2. Let $g\in{}^\o\GG_c$. Let $\JJ'\sneq\II$ be such that 
$\JJ\sub\JJ'$, $\o\in\Om_{\JJ'}$. We can find (and we fix) an increasing sequence $\GG_1\sub\GG_2\sub\do$ of
subsets of $\GG$ (depending on $\JJ'$) such that $\cup_{e\ge1}\GG_e=\GG$, $\GG_eP_{\JJ'}=\GG_e$ and 
$\GG_e/P_{\JJ'}$ is a projective variety in $\GG/P_{\JJ'}$ for each $e\ge1$. Let $e\ge1$. We have a 
commutative diagram of algebraic varieties with cartesian squares 
$$\CD \tM_{\JJ',g,e}@>\tj>>\tM_{\JJ',e}@<\tu<<E_e\T X@>\tq>>X\\
@Vp''VV     @Vp'VV      @V1\T pVV      @VpVV\\
M_{\JJ',g,e}@>j>>M_{\JJ',e}@<u<<E_e\T\bY@>q>>\bY\endCD$$
where 
$$X=\{(xP_\JJ,zU_{\JJ'})\in P_{\JJ'}/P_\JJ\T NP_{\JJ'}/U_{\JJ'};x\i zx\in NP_\JJ,p_\JJ(x\i zx)\in\fS\};$$
$\bY$ is the set of all $zU_{\JJ'}\in NP_{\JJ'}/U_{\JJ'}$ such that for some $xP_\JJ\in P_{\JJ'}/P_\JJ$, 
we have $x\i zx\in NP_\JJ$, $p_\JJ(x\i zx)\in cl(\fS)$;

$\tM_{\JJ',e}$ is the set of all pairs $(xP_\JJ,y(xU_{\JJ'}x\i))$ where $xP_\JJ\in\GG_e/P_\JJ$, 
$y(xU_{\JJ'}x\i)\in(xNP_\JJ x\i)/(xU_{\JJ'}x\i)$ are such that $p_\JJ(x\i yx)\in\fS$;

$M_{\JJ',e}$ is the set of all pairs $(xP_{\JJ'},y(xU_{\JJ'}x\i))$ where $xP_{\JJ'}\in\GG_e/P_{\JJ'}$, 
$y(xU_{\JJ'}x\i)\in(xNP_{\JJ'}x\i)/(xU_{\JJ'}x\i)$ are such that for some $v\in P_{\JJ'}$, we have
$y\in xvNP_\JJ v\i x\i$ and $p_\JJ(v\i x\i yxv)\in cl(\fS)$;
$$\tM_{\JJ',g,e}=\{(xP_\JJ,y(xU_{\JJ'}x\i))\in\tM_{\JJ',e};g\in xNP_\JJ x\i,g\i y\in xU_{\JJ'}x\i\};$$
$$M_{\JJ',g,e}=\{(xP_{\JJ'},y(xU_{\JJ'}x\i))\in M_{\JJ',e};g\in xNP_{\JJ'}x\i, g\i y\in xU_{\JJ'}x\i\};$$
$$E_e=\{hU_{\JJ'}\in\GG/U_{\JJ'};hP_{\JJ'}\in\GG_e/P_{\JJ'}\};$$

$j,\tj$ are the obvious imbeddings; $q,\tq$ are the obvious projections;
$$u(hU_{\JJ'},zU_{\JJ'})=(hP_{\JJ'},(hzh\i)(hU_{\JJ'}h\i)),$$ 
$$\tu(hU_{\JJ'},xP_\JJ,zU_{\JJ'})=(hxP_\JJ,(hzh\i)(hxU_{\JJ'}x\i h\i));$$
$$p(xP_\JJ,zU_{\JJ'})=zU_{\JJ'}, p'(xP_\JJ,y(xU_{\JJ'}x\i))=(xP_{\JJ'},y(xU_{\JJ'}x\i)),$$ 
$$p''(xP_\JJ,y(xU_{\JJ'}x\i))=(xP_{\JJ'},y(xU_{\JJ'}x\i)).$$
All the maps in the diagram are compatible with the natural actions of 
$\ov{P}_{\JJ'}$ where the action of
$\ov{P}_{\JJ'}$ on the four spaces on the left is trivial. Moreover, $\tu$ and $u$ are principal 
$\ov{P}_{\JJ'}$-bundles. Let 
$$\XX_{g,e}=\{xP_\JJ\in\GG_e/P_\JJ;x\i gx\in NP_\JJ,p_\JJ(x\i gx)\in\fS\}.$$
We define a local system $\hcs$ on $\XX_{g,e}$ by requiring that $\hcs_{xP_\JJ}=\cs_{p_\JJ(x\i gx)}$. We 
define a local system $\bcs$ on $X$ by requiring that $\bcs_{(xP_\JJ,zU_{\JJ'})}=\cs_{p_\JJ(x\i zx)}$. We 
define a local system $\dcs$ on $\tM_{\JJ',e}$ by requiring that 
$\dcs_{(xP_\JJ,y(xU_{\JJ'}x\i))}=\cs_{p_\JJ(x\i yx)}$. Note that $\hcs,\bcs,\dcs$ are well defined by the 
$\ov{P}_\JJ$-equivariance of $\cs$.

We have an isomorphism $\XX_{g,e}@>\si>>\tM_{\JJ',g,e}$ given by $xP_\JJ\m(xP_\JJ,g(xU_{\JJ'}x\i))$, under 
which these two varieties are identified; then the local system $\tj^*\dcs$ on $\tM_{\JJ',g,e}$ becomes 
$\hcs$. We have $\tq^*\bcs=\tu^*\dcs$ hence $q^*p_!\bcs=u^*(p'_!\dcs)$. The functors 
$$\cd_{\ov{P}_{\JJ'}}(\bY)@>q^*>>\cd_{\ov{P}_{\JJ'}}(E_e\T\bY)@<u^*<<\cd(M_{\JJ',e})@>j^*>>\cd(M_{\JJ',g,e})
$$
induce algebra homomorphisms
$$\align&\End_{\cd_{\ov{P}_{\JJ'}}(\bY)}(p_!\bcs)@>>>\End_{\cd_{\ov{P}_{\JJ'}}(E_e\T\bY)}(q^*p_!\bcs)@<<<
\End_{\cd(M_{\JJ',e})}(p'_!\dcs)\\&@>>>\End_{\cd(M_{\JJ',g,e})}(p''_!\tj^*\dcs)\endalign$$
of which the second one is an isomorphism since $u$ is a principal 
$\ov{P}_{\JJ'}$-bundle. Taking the 
composition of the first homomorphism with the inverse of the second one and with the third one and 
identifying $\End_{\cd_{\ov{P}_{\JJ'}}(\bY)}(p_!\bcs)=\End_{\cd(\bY)}(p_!\bcs)$ as in \cite{\CUSII, 1.16(a)}
we obtain an algebra homomorphism 
$$\End_{\cd(\bY)}(p_!\bcs)@>>>\End_{\cd(M_{\JJ',g,e})}(p''_!\tj^*\dcs).\tag a$$
For any $i\ge0$, $H^i_c(\XX_{g,e},\hcs)=H^i_c(\tM_{\JJ',g,e},\tj^*\dcs)=H^i_c(M_{\JJ',g,e},p''_!\tj^*\dcs)$ 
is naturally a module over the algebra $\End_{\cd(M_{\JJ',g,e})}(f_!\tj^*\dcs)$ hence, via (a), a module 
over the algebra $\End_{\cd(\bY)}(p_!\bcs)$. By 1.4(a), this last algebra may be canonically identified with
the group algebra $\bbq[{}^{\un{\JJ'-\JJ}}(\WW_{\JJ'})]$. Hence $H^i_c(\XX_{g,e},\hcs)$ is naturally a 
${}^{\un{\JJ'-\JJ}}(\WW_{\JJ'})$-module. Passing to the dual space we see that $H_{-i}(\XX_{g,e},\hcs)$ is 
naturally a ${}^{\un{\JJ'-\JJ}}(\WW_{\JJ'})$-module. Let 
$$\XX_g=\{xP_\JJ\in\GG/P_\JJ;x\i gx\in NP_\JJ,p_\JJ(x\i gx)\in\fS\}.$$ 
Now $\XX_g$ is the union of the increasing sequence of algebraic varieties $\XX_{g,1}\sub\XX_{g,2}\sub\do$ 
(the inclusions are imbeddings of algebraic varieties). Hence the notion of local system on $\XX_g$ is well 
defined. We shall denote again by $\hcs$ the local system on $\XX_g$ whose restriction to $\XX_{g,e}$ (for 
any $e\ge1$) is the local system denoted earlier by $\hcs$. The imbeddings above induce for any $i\ge0$ 
linear maps
$$H_{-i}(\XX_{g,1},\hcs)@>>>H_{-i}(\XX_{g,2},\hcs)@>>>H_{-i}(\XX_{g,3},\hcs)@>>>\do$$
which are compatible with the ${}^{\un{\JJ'-\JJ}}(\WW_{\JJ'})$-module structures. The direct limit of this 
system of linear maps is denoted by $H_{-i}(\XX_g,\hcs)$; it is a vector space independent of the choice of 
the sequence $\GG_1\sub\GG_2\sub\do$ and it carries a natural ${}^{\un{\JJ'-\JJ}}(\WW_{\JJ'})$-module 
structure. Moreover, if $g\in\GG_{rsc}$, this vector space is $0$ for $i$ large enough (depending on $g$) 
since

(b) {\it $\dim\XX_{g,e}$ is bounded as $e\to\iy$;}
\nl
a result closely related to (b) appears in \cite{\KL}, at least when $\JJ=\em$, but the general case can be 
reduced to the case where $\JJ=\em$.

We have a partition $\XX_g=\sqc_{\x\in\Om}{}^\x\XX_g$ where
$${}^\x\XX_g=\{xP_\JJ\in{}^\x\GG/P_\JJ;x\i gx\in NP_\JJ,p_\JJ(x\i gx)\in\fS\}.$$ 
Then for $i\ge0$, $H_{-i}({}^x\XX_g,\hcs)$ is defined in the same way as $H_{-i}(\XX_g,\hcs)$ and we have
canonically $H_{-i}(\XX_g,\hcs)=\op_{\x\in\Om}H_{-i}({}^\x\XX_g,\hcs)$.

\subhead 3.4\endsubhead
We preserve the setup of 3.2. From the classification of unipotent cuspidal character sheaves we see that
$\JJ$ has the following property:

(a) {\it for any $\o$-orbit $\KK$ (for the $\Om$-action on $\II$) such that $\KK\sneq\chJJ$, we have 
$\ss_\KK:=w_0^{\JJ\cup\KK}w_0^J=w_0^\JJ w_0^{\JJ\cup\KK}\in{}^{\un\KK}(\WW_{\JJ\cup\KK})$.}
\nl
From \cite{\CLA, 1.15} we see that the following holds.

(b) {\it If $\chJJ$ is a single $\o$-orbit then ${}^{\un{\chJJ}}\WW'=\{1\}$.}

(c) {\it If $\chJJ$ contains at least two $\o$-orbits then ${}^{\un{\chJJ}}\WW'$ is a Coxeter group (an 
irreducible affine Weyl group) on generators $\ss_\KK$ (one for each $\o$-orbit $\KK$ in $\chJJ$).}

Let $g\in{}^\o\GG_c$ and let $i\ge0$. We show:

(d) {\it There is a unique ${}^{\un{\chJJ}}\WW'$-module structure on $H_{-i}(\XX_g,\hcs)$ (see 3.3) such
that for any $\o$-orbit $\KK$ in $\chJJ$ with $\KK\ne\chJJ$, $\ss_\KK$ acts as in the
${}^{\un\KK}(\WW_{\JJ\cup\KK})$-module structure 3.3 on $H_{-i}(\XX_g,\hcs)$.}
\nl
If $\chJJ$ is a single $\o$-orbit, there is nothing to prove. Hence we can assume that $\chJJ$ contains at 
least two $\o$-orbits. It is enough to show that if $\KK\ne\KK'$ are distinct $\o$-orbits in $\chJJ$ and 
$\ss_\KK\ss_{\KK'}$ has finite order $m$ in ${}^{\un{\chJJ}}\WW'$ then the operators
$\ss_\KK,\ss_{\KK'}\in\Aut(H_{-i}(\XX_g,\hcs))$ defined by the ${}^{\un\KK}(\WW_{\JJ\cup\KK})$-module 
structure and by the ${}^{\un\KK'}(\WW_{\JJ\cup\KK'})$-module structure satisfy $(\ss_\KK\ss_{\KK'})^m=1$ as
operators in $\Aut(H_{-i}(\XX_g,\hcs))$. Since $m<\iy$ we see that $\KK\cup\KK'\sneq\chJJ$. Hence the 
${}^{\un{\KK\cup\KK'}}(\WW_{\JJ\cup\KK\cup\KK'})$-module structure on $H_{-i}(\XX_g,\hcs)$ is defined as 
in 3.3. From the definitions we see that this module structure restricts to the 
${}^{\un\KK}(\WW_{\JJ\cup\KK})$-module structure and to the ${}^{\un\KK'}(\WW_{\JJ\cup\KK'})$-module 
structure considered earlier. Hence it is enough to note that $(\ss_\KK\ss_{\KK'})^m=1$ as operators in 
$\Aut(H_{-i}(\XX_g,\hcs))$ in the ${}^{\un{\KK\cup\KK'}}(\WW_{\JJ\cup\KK\cup\KK'})$-module structure. This 
proves (d).

\subhead 3.5\endsubhead
We preserve the setup of 3.4. Let $\x\in\Om_\JJ$. We define $\t_\x:\XX_g@>>>\XX_g$ by $xP_\JJ\m xh\i P_\JJ$ 
where $h\in O_\x$; this is well defined and independent of the choice of $h$. The inverse image of the local
system $\hcs$ on $\XX_g$ is the local system $\hcs^\x$ defined in terms of $\cs^\x$ (see 3.2) in the same 
way as $\hcs$ is defined in terms of $\cs$. Using the isomorphism $\cs@>>>\cs^\x$ in 3.2, we can identify 
$\hcs^\x$ and $\hcs$. We see that $\t_\x$ induces for any $i\ge0$ an isomorphism 
$\t_{\x,*}:H_{-i}(\XX_g,\hcs)@>\si>>H_{-i}(\XX_g,\hcs)$. The operators $\t_{\x,*}$ for various 
$\x\in\Om_\JJ$ define a $\Om_\JJ$-module structure on $H_{-i}(\XX_g,\hcs)$. Combining this with the
${}^{\un{\chJJ}}\WW'$-module structure on $H_{-i}(\XX_g,\hcs)$ (see 3.4) we obtain a
${}^{\un{\chJJ}}\WW$-module structure on $H_{-i}(\XX_g,\hcs)$. (Note that ${}^{\un{\chJJ}}\WW$ is the 
semidirect product $\Om_J\cdot{}^{\un{\chJJ}}\WW'$ with ${}^{\un{\chJJ}}\WW'$ normal.)

From the definitions we see that for any $\x\in\Om$, $H_{-i}({}^\x\XX_g,\hcs)$ is a 
${}^{\un{\chJJ}}\WW'$-submodule of $H_{-i}(\XX_g,\hcs)$ and that the ${}^{\un{\chJJ}}\WW$-module 
$H_{-i}(\XX_g,\hcs)$ is induced by the ${}^{\un{\chJJ}}\WW'$-module $H_{-i}({}^1\XX_g,\hcs)$.

\subhead 3.6\endsubhead
We preserve the setup of 3.5. We now state:

(a) {\it Conjecture. For any $g\in{}^\o\GG_c$ and any $i\ge0$, the ${}^{\un{\chJJ}}\WW$-module 
$H_{-i}(\XX_g,\hcs)$ is finitely generated.}
\nl
If in addition, $g\in{}^\o\GG_{rsc}$ and $g$ is elliptic then $\XX_g$ is an ordinary algebraic variety (in 
the case where $\JJ=\em$ this follows from \cite{\KL, Cor.3.2}) hence $\dim H_{-i}(\XX_g,\hcs)<\iy$ and (a) 
holds. More generally, for any $g\in{}^\o\GG_{rsc}$ let $\L_g$ be the subgroup of $\GG$ consisting of all 
$x(\e)$ where $x$ runs over the 
one parameter subgroups $K^*@>>>\GG$ with image contained in the centralizer of $g$. Now $\L_g$ acts on 
$\XX_g$ by left multiplication and one can show that

(b) {\it $H_{-i}(\XX_g,\hcs)$ is finitely generated as a $\bbq[\L_g]$-module.}
\nl
(In the case where $\JJ=\em$ this follows from \cite{\KL, Prop.3.1(d)}.)

After a first version of this paper was posted, Zhiwei Yun pointed out to me 
that from \cite{\YU} one can deduce that (a) is a consequence of (b), at least if $\JJ=\em$, 
$g\in{}^\o\GG_{rsc}$ is topologically unipotent and the characteristic of $K$ is large enough.

\subhead 3.7\endsubhead
We preserve the setup of 3.5. Let $\Om_{\JJ,1}$ (resp. $\Om_{\JJ,2}$) be the kernel (resp. image) of the 
homomorphism $f:\Om_\JJ@>>>\Aut({}^{\un{\chJJ}}\WW')$, $\x\m[w\m\x w\x\i]$. We can form the semidirect 
product $\Om_{\JJ,2}\cdot{}^{\un{\chJJ}}\WW'$ using the action of $\Om_{\JJ,2}$ on ${}^{\un{\chJJ}}\WW'$
defined by $f$. From the classification of arithmetic diagrams in \cite{\CLA, \S7} we see that we have

(i) $\Om_{\JJ,1}=\{1\},\Om_{\JJ,2}=\Om_J$, or

(ii) $\Om_{\JJ,1}=\Om_J,\Om_{\JJ,2}=\{1\}$, or

(iii) $\Om_{\JJ,1}\cong\Om_{\JJ,2}\cong\ZZ/2\ZZ$, $\Om_\JJ\cong(\ZZ/2\ZZ)\op(\ZZ/2\ZZ)$.
\nl
(Case (iii) occurs only in the case 7.44 and 7.45 in \cite{\CLA, \S7} with $\GG$ of type $D_n$, $n$ even,
$n>4$. In this case $\Om_{\JJ,1}=\{1,\o\}$.)
We define an isomorphism $\Om_{\JJ,1}\T\Om_{\JJ,2}@>\si>>\Om_\JJ$ as follows. In case (i) or (ii) this is the
obvious isomorphism. In case (iii) the isomorphism restricted to $\Om_{\JJ,1}$ is the obvious imbedding, 
while the image of the nontrivial element of $\Om_{\JJ,2}$ is the unique element $\g\in\Om_\JJ$ such that 
the permutation $i\m\g(i)$ of $\II$ has at least two fixed points. Using this isomorphism, the group 
${}^{\un{\chJJ}}\WW=\Om_J\cdot{}^{\un{\chJJ}}\WW'$ becomes the direct product of groups 
$\Om_{\JJ,1}\T(\Om_{\JJ,2}\cdot{}^{\un{\chJJ}}\WW')$. Hence the algebra $\bbq[{}^{\un{\chJJ}}\WW]$ becomes 
the algebra $\bbq[\Om_{\JJ,1}]\ot\bbq[\Om_{\JJ,2}\cdot{}^{\un{\chJJ}}\WW']$. Now the group algebra 
$\bbq[\Om_{\JJ,1}]$ is canonically isomorphic to the direct sum of copies of the algebra $\bbq$ (indexed by 
the characters $\Om_{\JJ,1}@>>>\bbq^*$). Hence

(a) {\it the algebra $\bbq[{}^{\un{\chJJ}}\WW]$ is canonically isomorphic to the direct sum of copies
of the algebra $\bbq[\Om_{\JJ,2}\cdot{}^{\un{\chJJ}}\WW']$ (indexed by the characters 
$\Om_{\JJ,1}@>>>\bbq^*$).}

\subhead 3.8\endsubhead
Let $\cg$ be a connected simply connected almost simple group over $\CC$ of type dual to that of $G$. We 
have a natural bijection
$$\io:\Om@>\si>>\Hom(\cz_\cg,\CC^*).\tag a$$
For any $\o\in\Om$ we set $\fZ^\o=\fZ_{\cg,\io(\o)}$, see 2.2.

Let $\z=(\fc,\fF)\in\fZ_\cg$. Now $\z$ corresponds under 2.6(c) to a quadruple $(J,\cc,\ce,E)$ with
$E=E_\z\in\Irr\tW_{\chJ}$ (notation of 2.7). By \cite{\CLA, 6.3}, \cite{\CLA,\S7} to the triple 
$(J,\cc,\ce)$ corresponds a triple $(\JJ,\o,A)\in\fV$ so that $\o$ satisfies $\z\in\fZ^\o$ and 
$\Om_{\JJ,2}\cdot{}^{\un{\chJJ}}\WW'$ is identified with the group $\tW_{\chJ}$. 
Thus $E_\z$ and $\hE_\z$ (see 2.7) become representations of one specific copy
of $\bbq[\Om_{\JJ,2}\cdot{}^{\un{\chJJ}}\WW']$ in the direct sum of algebras in 3.7(a); hence they can be
viewed as representations of $\bbq[{}^{\un{\chJJ}}\WW]$ on which the other copies of 
$\bbq[\Om_{\JJ,2}\cdot{}^{\un{\chJJ}}\WW']$ in the direct sum in 3.7(a) act as zero. These representations
of $\bbq[{}^{\un{\chJJ}}\WW]$ are denoted by $\uE_\z$ and $\huE_\z$. We now define for any 
$g\in{}^\o\GG_c$ and any $i\ge0$ a $\bbq$-vector space 
$$\AA^{\z;i}_g=\Hom_{{}^{\un{\chJJ}}\WW}(H_{-i}(\XX_g,\hcs),\huE_\z)$$
(we regard $H_{-i}(\XX_g,\hcs)$ as a ${}^{\un{\chJJ}}\WW$-module as in 3.5). Note that for any fixed 
$g\in{}^\o\GG_{rsc}$, 

(b) {\it $\AA^{\z;i}_g$ is zero for all but finitely many $i$.}
\nl
(This follows from 3.3(b).) 

Now let $g\in{}^\o\GG_c$ and let $h\in\GG$. Define $\nu_h:\XX_g@>\si>>\XX_{hgh\i}$ by 
$xP_\JJ\m hxP_\JJ$. If $\hcs'$ is the local system on $\XX_{hgh\i}$ defined in the same way as $\hcs$ on 
$\XX_g$ then $\nu_h^*\hcs'=\hcs$. Hence $\nu_h$ induces for any $i\ge0$ an isomorphism 
$[\nu_h]:H_{-i}(\XX_g,\hcs)@>\si>>H_{-i}(\XX_{hgh\i},\hcs')$. From the definitions we see that $[\nu_h]$ is 
compatible with the ${}^{\un{\chJJ}}\WW$-module structures. Hence its inverse induces a linear isomorphism 

$\AA^{\z;i}_g@>\si>>\AA^{\z;i}_{hgh\i}$ 
\nl
for any $\z\in\ufZ_\cg$. Thus the collection of vector spaces $\AA^{\z;i}_g$ (with $g\in{}^\o\GG_c$) has a 
natural $\GG$-equivariant structure.

We conjecture:

(e) {\it for fixed $\z$ and $i$, there exists an $l$-adic constructible sheaf $\cf$ on ${}^\o\cx$ (see 3.1)
such that the stalk of $\cf$ at $(g,B)\in {}^\o\cx$ is the vector space $\AA^{\z;i}_g$.}
\nl
This conjecture implies the inequality 

(f) {\it $\dim\AA^{\z;i}_g<\iy$ for any $g\in{}^\o\GG_c$.}
\nl
This would also follow from Conjecture 3.6(a). 

\subhead 3.9\endsubhead
Until the end of 3.12 we assume that $\kk$ is an algebraic closure of the finite field $\FF_q$.
Let $K_0=\FF_q((\e))$, a subfield of $K$. Let $F:\GG@>>>\GG$ be the bijective homomorphism induced by the 
map $x\m x^q$ of $\kk$ into itself so that the fixed point set $\GG^F$ is equal to $G(K_0)$ for a split 
$K_0$-form of our group. 

Let $\o\in\Om$ and let $\z=(\fc,\fF)\in\fZ^\o$. Our next objective is to associate to $\z$ a function 
$\tt_\z:{}^\o\GG_{rsc}^F@>>>\bbq\cup\{\iy\}$. (Here ${}^\o\GG_{rsc}^F:={}^\o\GG_{rsc}\cap\GG^F$.) Define 
$(\JJ,\o,A)\in\fV$ in terms of $\z$ as in 3.8. 
Define $\fS,\cs$ in terms of $A$ as in 3.2. We have $F(P_\JJ)=P_{\JJ}$ and $F$ induces a 
morphism of algebraic groups $\ov{NP}_\JJ@>>>\ov{NP}_\JJ$ denoted again by $F$; this is the Frobenius map 
for an $\FF_q$-rational structure on $\ov{NP}_\JJ$ such that every connected component of $\ov{NP}_\JJ$ is 
defined over $\FF_q$. We have necessarily $F(\fS)=\fS$ and there exists an isomorphism 
$\ps:F^*\cs@>\si>>\cs$ of local systems over $\fS$. We shall assume (as we may) that for any $\s\in\fS$ and 
any $m\ge1$ such that $F^m(\s)=\s$, the eigenvalues of $F^m$ on $\cs_\s$ are roots of $1$; moreover if 
$\JJ=\em$ so that $\cs=\bbq$ we shall take $\ps=1$.

For any $g\in{}^\o\GG_c$ let $\XX_g$ be as in 3.3. We have $F(g)\in{}^\o\GG_c$ and we define 
$\ph:\XX_g@>>>\XX_{F(g)}$ by $xP_\JJ\m F(x)P_\JJ$. Let $\hcs'$ be the local system on $\XX_{F(g)}$ defined 
in the same way as $\hcs$ on $\XX_g$, see 3.3. Then $\ps$ induces an isomorphism $\hcs@>\si>>\ph^*\hcs'$. 
Hence for any $i\ge0$ there is an induced isomorphism 
$$H_{-i}(\XX_g,\hcs)@>\si>>H_{-i}(\XX_{F(g)},\hcs').\tag a$$
From the definitions we see that this isomorphism is compatible with the ${}^{\un{\chJJ}}\WW$-module 
structures. Now (a) induces for any $g\in{}^\o\GG_c$ and any $i\ge0$ an isomorphism 
$\ti\ps_g:\AA^{\z;i}_{F(g)}@>>>\AA^{\z;i}_g$. If $F(g)=g$ then $\tr(\ti\ps_g:\AA^{\z;i}_g@>>>\AA^{\z;i}_g)$ 
is well defined (assuming that 3.8(f) holds); we denote it by $\tt_{\z;i}(g)$. 

Assuming now that $g=F(g)\in{}^\o\GG_{rsc}$, we set 
$\tt_\z(g)=\sum_{i\ge0}(-1)^i\tt_{\z;i}(g)$. Thus $g\m\tt_\fz(g)$ is a well defined function on 
${}^\o\GG_{rsc}^F$ with values in $\bbq$, except that at points $g$ where 3.8(f) fails for some $i$ we set 
$\tt_\z(g)=\iy$. The function $\tt_\z:{}^\o\GG_{rsc}^F@>>>\bbq\cup\{\iy\}$ is constant on each 
$\GG^F$-conjugacy class in ${}^\o\GG_{rsc}^F$ (this follows from the compatibility of $\ti\ps$ above with 
the $\GG$-equivariant structure on $\AA^{\z;i}$). In the remainder of this section we assume that 3.8(f) 
holds so that $\tt_\z:{}^\o\GG_{rsc}^F@>>>\bbq$. We shall also regard $\tt_\z$ as a function 
$\GG_{rsc}^F@>>>\bbq$, equal to zero on $\GG_{rsc}^F-{}^\o\GG_{rsc}^F$. (We set 
$\GG_{rsc}^F=\GG_{rsc}\cap\GG^F$.) The functions $\tt_\z$ on $\GG_{rsc}^F$ (for $\z\in\fZ^\o$, $\o\in\Om$ 
variable) are said to be the {\it unipotent almost characters} of $\GG^F$. They are defined up to 
multiplication by a root of $1$. Let $\VV''$ be the subspace of the vector space of class functions 
$\GG_{rsc}^F@>>>\bbq$ generated by the unipotent almost characters of $\GG^F$.

\subhead 3.10\endsubhead
Let $\cu$ be the set of isomorphism classes of unipotent representations of $\GG^F$ (see \cite{\CLA, 0.3}).
Recall that \cite{\CLA, 6.5} gives a bijection 
$$\fZ^1\lra\cu\tag a$$
with $\fZ^1$ as in 3.8. Now for $\z\in\fZ^1$ the irreducible $\GG^F$-module $R_\z$ corresponding to
$\z$ under (a) is canonically the 
quotient of a $\GG^F$-module $\hR_\z$ which is standard in the sense of Langlands and which has finite
length with all composition factors being again unipotent representations. Hence $\hR_\z$ has a well defined
character (in the sense of Harish-Chandra) whose restriction to $\GG_{rsc}^F$ is denoted by $\ph_\z$. We 
call the class functions $\ph_\z:\GG_{rsc}^F@>>>\bbq$ (for $\z\in\fZ^1$) the {\it unipotent characters} of 
$\GG^F$. Let $\VV'$ be the subspace of the vector space of class functions $\GG_{rsc}^F@>>>\bbq$ generated 
by the unipotent characters.

(a) {\it Conjecture. We have $\VV'=\VV''$.}

\subhead 3.11\endsubhead
We now formulate a refinement of conjecture 3.10(a). We fix $\o\in\Om$ and let $\uc$ be a unipotent class of
$\cg$. Let ${}^{\uc}\fZ^\o$ be the set of all $(\fc,\fF)\in\fZ^\o$ such that the unipotent part of any 
element of $\fc$ lies in $\uc$. 

Let $\VV'_{\o,\uc}$ be the subspace of the vector space of class functions ${}^\o\GG_{rsc}^F@>>>\bbq$ 
generated by the restrictions $\ph_\z^\o:=\ph_\z|_{{}^\o\GG_{rsc}^F}$ (for various $\z\in{}^{\uc}\fZ^1$). Let
$\VV''_{\o,\uc}$ be the subspace of the vector space of class functions ${}^\o\GG_{rsc}^F@>>>\bbq$ generated
by the unipotent almost characters $\tt_\z$ (for various $\z\in{}^{\uc}\fZ^\o$).

(a) {\it Conjecture. We have $\VV'_{\o,\uc}=\VV''_{\o,\uc}$.}

\subhead 3.12\endsubhead
Let us now replace the Frobenius map $F:\GG@>>>\GG$ by another Frobenius map $F':\GG@>>>\GG$ such that the 
fixed point set $\GG^{F'}$ is an inner form of the split form of $G$ over $K_0$. Then the analogue of $\cu$ 
for $\GG^{F'}$ is well defined (it is in bijection with $\fZ^\x$ for some $\x\in\Om$ which depends on $F'$, 
see \cite{\CLA}). The unipotent almost characters on $\GG^{F'}$ are again defined by taking trace of 
$F'$ instead of $F$ on the same geometric objects which were used for $\GG^F$. We expect that the analogues 
of 3.10(a) and 3.11(a) continue to hold.

\head 4. Examples\endhead
\subhead 4.1\endsubhead
In this section we preserve the notation of \S3. Let $\co=\kk[\e]]$, $\co^*=\co-\e\kk[[\e]]$. Let 
$v:K^*@>>>\ZZ$ be the group homomorphism such that $v(\e^n)=n$ for $n\in\ZZ$ and $v(\co^*)=0$. We set
$v(0)=\iy$.

\subhead 4.2\endsubhead
Let $\GG_{cvr}$ be the set of compact very regular elements in $\GG$ (see \cite{\KmL}). We have 
$\GG_{cvr}\sub{}^1\GG_{rsc}$. Let $g\in\GG_{cvr}$. We show that 3.6(a) holds for this $g$.
If $(\JJ,\o,A)\in\fV$ is such that $\JJ\ne\em$, $\o=1$ then $\XX_g=\em$. Now assume that 

(a) $(\JJ,\o,A)\in\fV$ is such that $\JJ=\em$, $\o=1$ (hence $A$ is $\bbq$ up to a shift).
\nl
From \cite{\KmL, 2.2} we see that 

(b) {\it $H_{-i}(\XX_g,\bbq)=0$ if $i>0$, $H_0(\XX_g,\bbq)$ is a vector space with basis $\b$ and the module 
structure over ${}^{\un{\chJJ}}\WW=\WW$ is such that $\WW$ permutes $\b$ simply transitively.}
\nl
In particular, 3.6(a) holds for this $g$.

Now let $\z\in\fZ^1$ be such that the corresponding triple $(\JJ,\o,A)$ (see 3.8) is as in (a). Using (b) we
see that if $g\in\GG_{cvr}$ then $\AA^{\z;i}_g=0$ if $i>0$; moreover any $b\in\b$ (see (b)) defines an 
isomorphism of vector spaces $\AA^{\z;0}_g@>\si>>\huE_\z$ which carries any element of $\AA^{\z;0}_g$ 
(viewed as a linear map $H_0(\XX_g,\bbq)@>>>\huE_\z$) to the image of $b$ under this linear map.

From this we see that if we are in the setup of 3.9 and $g$ is in addition fixed by $F$ (so that $g$ 
corresponds as in \cite{\KmL} to an element of finite order $w\in\WW'$ defined up to conjugacy) then

(c) $\tt_\z(g)=\tr(w,\huE_\z)$.

\subhead 4.3\endsubhead
Now assume that $G=PGL_2$ and 

(a) $g\in\GG$ is represented by $\left(\sm a&b\\c&d\esm\right)$ where $a,b,d\in1+\e\co$, $c\in\e+\e^2\co$.
\nl
Let $\z\in\fZ^1$. Then the corresponding triple $(\JJ,\o,A)$ (see 3.8) is as in 4.2(a). In this case we have
$H_{-i}(\XX_g,\bbq)=0$ if $i>0$ and $H_0(\XX_g,\bbq)$ is a vector space with basis $\{b_1,b_2\}$ and the 
module structure over ${}^{\un{\chJJ}}\WW=\WW$ is such that $\WW'$ acts trivially and any element in 
$\WW-\WW'$ interchanges $b_1,b_2$. It follows that $\AA^{\z;i}_g=0$ if $i>0$ and $\AA^{\z;0}_g$ is 
isomorphic to the space of $\WW'$-invariant elements in $\huE_\z$.

\subhead 4.4\endsubhead
Now assume that $G=PGL_2$ and 

(b) $g\in\GG$ is represented by $\left(\sm1+\e&0\\0&1-\e\esm\right)$.
\nl
Let $\z\in\fZ^1$. Then the corresponding triple $(\JJ,\o,A)$ (see 3.8) is as in 4.2(a). In this case 
$H_0({}^1\XX_g,\bbq)$ is one-dimensional and has trivial action of $\WW'$. For $i>0,i\ne2$ we have
$H_{-i}({}^1\XX_g,\bbq)=0$. Now $H_{-2}({}^1\XX_g,\bbq)$ has basis $\{b_n;n\in\ZZ\}$ and the simple 
reflection $s_1,s_2$ of $\WW'$ acts as follows.
$$s_ib_n=-b_n\text{ if }n=i\mod2,\qua s_ib_n=b_n+b_{n-1}+b_{n+1}\text{ if }n=i+1\mod2.$$
We see that $b_0,b_1$ generate the $\WW'$-module $H_{-2}({}^1\XX_g,\bbq)$. In particular, 3.6(a) holds for 
our $g$.

It follows that $\AA^{\z;i}_g=0$ if $i>0,i\ne2$, $\AA^{\z;0}_g$ is isomorphic to the space of 
$\WW'$-invariant elements in $\huE_\z$ and $\AA^{\z;2}_g$ is isomorphic to the space of sequences 
$(x_n)_{n\in\ZZ}$ with $x_n\in\huE_\z$ satisfying

$s_ix_n=-x_n$ if $n=i\mod2$, $s_ix_n=x_n+x_{n-1}+x_{n+1}$ if $n=i+1\mod2$;
\nl
here $s_i (i=1,2)$ are the simple reflections of $\WW'$.

We now assume that $\z=(\{1\},\CC)$. In this case $\huE_\z$ has a basis $\{x,x'\}$ such that the action of
$s_i$ is as follows:

$s_ix=-x$ for $i=1,2$, $s_1x'=x'+x$, $s_2x'=x'$. 
\nl
Thus the sign representation of $\WW'$ is a submodule and the unit representation of $\WW'$ is a quotient
module. We see that $\AA^{\z;i}_g=0$ for all $i\ne2$. Moreover $\AA^{\z;2}_{g'}$ can be identified (via 
$(x_n)\m(u_n), x_n=u_nx$) with the vector space of sequences $(u_n)_{n\in\ZZ}$ with $u_n\in\CC$ satisfying 
$-u_n=u_n+u_{n-1}+u_{n+1}$ for $n\in\ZZ$. This is a two-dimensional vector space: if $u_0,u_1$ are given 
then the other $u_n$ are uniquely determined. Thus $\dim\AA^{\z;2}_g=2$. From these computations we see also
that if we are in the setup of 3.9 then

$\tt_\z(g)=2q$.
\nl
Since the value of the Steinberg character of $\GG^F$ at $g$ is $2q-1$ (see \cite{\Kml}) we see that 
$\ph_\z^1(g)=2q$ (note that $2q=(2q-1)+1$ where the last $1$ comes from the unit representation) 
so that $\ph_\z^1(g)=\tt_\z(g)$. (We expect that this equality holds for any 
$g\in{}^1\GG_{rsc}^F$.)

Note that in the arguments above, our choice of a co-standard module plays a key role; if we replace it by 
one in which the unit representation of $\WW'$ is a submodule and the sign representation of $\WW'$ is a 
quotient module, the arguments above would collapse.

\subhead 4.5\endsubhead
In this subsection we assume that $G=PGL_2$ and $\o\in\Om$ is not $1$; we shall prove that 
3.11(a) holds in this case. Let $\tGG=GL_2(K)$ and let $\p:\tGG@>>>\GG$ be the obvious homomorphism. We set
$$\tGG^1=\{\left(\sm a&b\\c&d\esm\right)\in\tGG;v(ad-bc)=\text{even}\},
\tGG^2=\{\left(\sm a&b\\c&d\esm\right)\in\tGG;v(ad-bc)=\text{odd}\},$$
$$I^1=\{\left(\sm a&b\\c&d\esm\right)\in\tGG;v(a)=v(d)=m,v(b)\ge m,v(c)>m \text{ for some }m\in\ZZ\},$$
$$I^2=\{\left(\sm c&d\\a&b\esm\right);v(a)=v(d)+1=m+1,v(b)\ge m+1,v(c)\ge m+1\text{ for some }m\in\ZZ\}.$$
Then $\tGG=\tGG^1\sqc\tGG^2$, $I^1\sub\tGG^1$, $I^2\sub\tGG^2$
 and $\p(I^1)$ is an Iwahori subgroup of $\GG$ hence we can assume that 
$P_\em=\p(I_1)$. We have $NP_\em=\p(I_1\cup I_2)$ and ${}^\o NP_\em=\p(I_2)$. Let $\un B$ be the image under
$\p$ of the group of upper triangular matrices in $\tGG$. We show:

(a) {\it If $g\in{}^\o\GG$ normalizes some Iwahori subgroup
then $g$ is not contained in any $\GG$-conjugate of $\un B$; in particular we have $g\in\GG_{rsc}$.}
\nl
We can assume that $g\in\p(I^2)$ that is $g=\p(\tg)$, where 
$\tg=\left(\sm \e c&d\\ \e a&\e b\esm\right)$ and $d\in\co^*$, $a\in\co^*$, $b\in\co$, $c\in\co$. Assume 
that $g$ is contained in a $\GG$-conjugate of $\un B$; then the eigenvalues $\l,\l'$ of $\tg$ are contained 
in $K$. We have $\l+\l'=\e(b+c)$, $\l\l'=\e^2bc-\e ad$ hence 
$(\l-\l')^2=(\l+\l')^2-4\l\l'=\e^2(b+c)^2-4\e^2bc-\e ad$ so that $v((\l-\l')^2)=1$ contradicting the 
inclusion $\l-\l'\in K$. This proves (a).

We show:

(b) {\it If $\tg\in I^2$, $\tg'\in I^2$, $h\in\tGG^1$ satisfy $\tg'=h\i\tg h$ then $h\in I^1$.}
\nl 
We can assume that 
$$\tg=\left(\sm\e c&d\\ \e a&\e b\esm\right),\qua \tg'=\e^m\left(\sm\e c'&d'\\ \e a'&\e b'\esm\right),$$
where $a,d,a',d'$ are in $\co^*$ and $b,c,b',c'$ are in $\co$. We have $\det(\tg)=\det(\tg')$ where 
$v(\det(\tg))=1$, $v(\det(\tg'))=2m+1$. Hence $m=0$. Since $\tg,\tg'$ are conjugate we have $b+c=b'+c'$, 
$ad-\e bc=a'd'-\e b'c'$. Let $r=\left(\sm 1&(c'-c)/a\\ 0&a'/a\esm\right)$. We have $r\in I^1$, $r\tg=\tg'r$.
Thus $\tg=(hr)\i\tg hr$. If we can prove that $hr\in I^1$ it would follow that $h\in I^1$. Thus we are 
reduced to the case where $\tg'=\tg$. We write 
$$h=\left(\sm x&y\\z&u\esm\right)$$
with $x,y,z,u\in K$, $v(xu-yz)\in2\ZZ$. From $h\tg=\tg h$ we have $ay\e=dz$, $xd+yb\e=cy\e+du$ hence 
$z=\e ya/d$, $u=x+\e y(b-c)/d$. Thus
$\det(h)=x^2+\e xy(b-c)/d-\e y^2 a/d$. If $x=0$ then $y\ne0$ and $\det h=-\e y^2 a/d$, so that
$v(\det h)$ is odd, contradicting $h\in\tGG^1$. Thus $x\ne0$ so that $M:=v(x)<\iy$. If $v(y)<M$ then
$v(x^2)=2M$, $v(\e xy(b-c)/d)<M+1+v(y)\le2M$, $v(\e y^2 a/d)=2v(y)+1<2M$. We see that $v(\det h)=2v(y)+1$ is
odd, contradicting $h\in\tGG^1$. Thus we have $v(y)\ge M$. We have $v(x)=M$, $v(y)\ge M$, 
$v(\e ya/d)\ge M+1$, $v(x+\e y(b-c)/d)=v(x)=M$. Thus $h\in I^1$. This proves (b).

We now show: 

(c) {\it Let $g\in{}^\o\GG^2$ be such that $g$ normalizes some Iwahori subgroup. Then the set 
$\{xP_\em\in\GG/P_\em;x\i gx\in{}^\o NP_\em\}$ consists of two points.}
\nl
An equivalent statement is as follows. Let $g\in\tGG^2$ be such that some $\tGG$-conjugate of $g$ is in 
$I^2$. Then the set $\{xI^1\in\tGG/I^1;x\i gx\in I^2\}$ consists of two points.

We can assume that $g\in I^2$. Then our set contains at least two points, $I^1$ and $I^2$. If $x\in\tGG^1$, 
$x\i gx\in I^2$ then by (b) we have $x\in I^1$ so that $xI^1=I^1$. If $x\in\tGG^2$, $x\i gx\in I^2$ then, 
choosing $s\in I^2$, we have $(xs)\i g(xs)\in I^2$, $xs\in\tGG^1$ hence using again (b) we have $xs\in I^1$ 
so that $xI^1=I^2$. This proves (c).

For $g$ as in (c) the set $\XX_g$ in 3.3 (with $\JJ$ necessarily empty) has exactly two elements. Hence 
3.6(a) holds in this case; hence 3.8(f) also holds. 

Assume that we are in the setup of 3.9.
If $\z=(\fc,\fF)\in\fZ^1$ is such that $\fc$ is a semisimple class of $\cg$ then $\hR_\z$ in 3.10 is induced
from a character of $\un B\cap\GG^F$; using (a) we see that $\ph_\z^\o$ is identically zero on 
${}^\o\GG_{rsc}^F$. If $\z=(\fc,\fF)\in\fZ^1$ is such that $\fc$ is the regular unipotent class of $\cg$ (so
that $\fF=\bbq$) then $\hR_\z$ in 3.10 is the unit representation of $\GG^F$; it follows that $\ph_\z^\o$ is
identically $1$ on ${}^\o\GG_{rsc}^F$. If $\z=(\fc,\fF)\in\fZ^1$ is such that $\fc$ is $-1$ times the 
regular unipotent class of $\cg$ (so that $\fF=\bbq$) then $\hR_\z$ in 3.10 is the one dimensional 
representation of $\GG^F$ which is trivial on ${}^1\GG^F$ and on which elements in ${}^\o\GG^F$ acts as 
$-1$; it follows that $\ph_\z^\o$ is identically $-1$ on ${}^\o\GG_{rsc}^F$.
We see that $\VV'_{\o,\uc}$ in 3.11 is $0$ if $\uc=\{1\}$ and is the vector space spanned by the constant 
function $1$ if $\uc$ is the regular unipotent class of $\cg$.

If $\z=(\fc,\fF)\in\fZ^\o$ then $\fc$ is the regular unipotent class of $\cg$ times $1$ or $-1$ and $\fF$ is
a cuspidal local system on $\fc$. The corresponding $(\JJ,\o,A)\in\fV$ satisfies $\JJ=\em$, we have 
${}^{\un{\chJJ}}\WW=\bbq[\Om]$ and $\huE_\z$ is a one dimensional representation of $\Om$. Using (c) we see 
that for any $g\in{}^o\GG_{rsc}$ we have $\AA^{\z;0}_g=\bbq$, $\AA^{\z;i}_g=0$ if $i>0$. It also follows 
that $\tt_\z$ is a constant (equal to a root of $1$) on ${}^\o\GG_{rsc}^F$. We see that $\VV''_{\o,\uc}$ in 
3.11 is $0$ if $\uc=\{1\}$ and is the vector space spanned by the constant function $1$ if $\uc$ is the 
regular unipotent class of $\cg$. Thus 3.11(a) holds in this case.

\subhead 4.6\endsubhead
We return to the general case. 
In this subsection we take $\o=1$ and $\uc$ to be the regular unipotent class in $\cg$. We show that 3.11(a)
suggests a geometric property of $\GG$. Now if $\z\in{}^{\uc}\fZ^1$ then $\ph_\z^1$ is identically $1$; 
hence $\VV'_{1,\uc}$ is spanned by the constant function $1$ on ${}^1\GG_{rsc}^F$. Assuming 3.11(a), we 
deduce that if $\z=(\uc,\CC)$ then $\tt_\z$ is a constant function on ${}^1\GG_{rsc}^F$. This suggests that
for any $g\in{}^1\GG_{rsc}$ we have

(a) {\it $\AA^{\z,i}_g=0$ if $i>0$ and $\AA^{\z,0}_g\cong\bbq$.}
\nl
In this case $\AA^{\z,i}$ is defined in terms of a triple as in 4.2(a). Hence (a) implies:

(b) {\it for any $g\in{}^1\GG_{rsc}$, the space of $\WW$-coinvariants in $H_{-i}(\XX_g,\bbq)$ is $0$ if 
$i>0$ and is one dimensional if $i=0$.}
\nl
Note that in the setup of 4.4 this follows directly from the results in 4.4.

\subhead 4.7\endsubhead
Assume that we are in the setup of 3.9. Let $\uc$ be a unipotent class in $\cg$.
Let $u\in\uc$ and let $\G$ be a maximal reductive subgroup of $Z_\cg(u)$. 
Let $M(\G)$ be the set consisting of all pairs $(x,\s)$ where $x\in\G$ is semisimple and is defined up to 
$\G$-conjugacy and $\s\in\Irr Z_\G(x)/Z_\G(x)^0$. Let $M(\G)^1$ be the set of all $(x,\s)\in M(\G)$ such 
that the restriction of $\s$ to $\cz_\cg$ is trivial. We have a bijection

(a) ${}^{\uc}\fZ^1\lra M(\G)^1$;
\nl
to $(\fc,\fF)\in{}^{\uc}\fZ^1$ corresponds to $(x,\s)\in M(\G)^1$ if $ux\in\fc$ and 
$Z_\cg(ux)/Z_\cg(ux)^0=Z_\G(x)/Z_\G(x)^0$ acts on the stalk $\fF_{ux}$ through a representation isomorphic 
to $\s$. For $(x,\s)\in M(\G)^1$ we write $\tt_{x,\s}$, $\hR_{x,\s}$, $\ph_{x,\s}$ instead of $\tt_\z$, 
$\hR_\z$, $\ph_\z^1$ (see 3.9, 3.10, 3.11) where $\z\in{}^{\uc}\fZ^1$ corresponds as above to $(x,\s)$. 

We expect that if $(x,\s)\in M(\G)^1$ is such that $Z_\G(x)/\cz_{\cg}$ is connected (so that $\s=1$) then 

(b) $\tt_{x,\s}=\ph_{x,\s}$ up to a root of $1$.
\nl
In general, we expect that the functions $\tt_{x,\s}$, $\ph_{x,\s}$ on ${}^1\GG_{rsc}^F$ are related by a 
version of the nonabelian Fourier transform \cite{\FOUR}. 
To illustrate this we consider the case where $\cg$ is of type $B_2$ or $G_2$
and $\uc$ is the subregular unipotent class in $\cg$. We can identify $\G=\CC^*\cdot<r>$, $\G=S_3$ if $\cg$ 
is of type $B_2,G_2$ respectively. Here $<r>$ is a group of order $2$ with generator $r$ acting on $\CC^*$ 
by $z\m z\i$ and $S_3$ is the symmetric group in three letters. We expect that for any $(x,\s)\in M(\G)^1$ 
the identities (c)-(d) hold (up to multiplication by a root of $1$) for the functions 
$\tt_{x,\s}$, $\ph_{x,\s}$ on ${}^1\GG_{rsc}^F$. 

(c) $\tt_{x,\s}=\ph_{x,\s}$
\nl
if $\cg$ is of type $B_2$ and $x\in\G^0,x^2\ne1$ (this is a special case of (b));

(d) $\tt_{1,1}=(\ph_{1,1}+\ph_{r,1}+\ph_{1,\e}+\ph_{r,\e})/2$,

$\tt_{r,1}=(\ph_{1,1}+\ph_{r,1}-\ph_{1,\e}-\ph_{r,\e})/2$,

$\tt_{1,\e}=(\ph_{1,1}-\ph_{r,1}+\ph_{1,\e}-\ph_{r,\e})/2$,

$\tt_{r,\e}=(\ph_{1,1}-\ph_{r,1}-\ph_{1,\e}+\ph_{r,\e})/2$,
\nl
if $\cg$ is of type $B_2$ and $\e$ is the nontrivial character of $\G/\G^0$;

(e) $\tt_{x,\s}=\sum_{(x',\s')\in M(\G)}\{(x,\s),(x',\s')\}\ph_{x',\s'}$
\nl
if $\cg$ is of type $G_2$ (here $\{(x,\s),(x',\s')\}$ is the nonabelian Fourier transform matrix 
\cite{\FOUR}). Note that in the setup of (d) we have

(f) $\ph_{-1,1}=\ph_{1,1}$, $\ph_{-1,\e}=\ph_{1,\e}$, 

$\tt_{-1,1}=\tt_{1,1}$, $\tt_{-1,\e}=\tt_{1,\e}$.
\nl
Clearly, if (c)-(e) hold, then 3.11(a) holds in our case. Now we can check that (c)-(e) hold after 
restriction to $\GG_{cvr}^F=\GG_{cvr}\cap\GG^F$, using 4.2(c), the known results about the character of 
unipotent representations of reductive groups over $\FF_q$ and the results of \cite{\KmL}. (See in
particular \cite{\KmL, 5.3}.) This provides additional support for 3.11(a).

\head 5. Affine flag manifolds in unequal characteristic \endhead
\subhead 5.1\endsubhead
Let $\kk$ be an algebraic closure of the finite field $\FF_q$ (of characteristic $p$), let $\co$ be the ring
of Witt vectors over $\kk$ and let $K$ be the quotient field of $\co$. Let $\GG$ be the group of 
$K$-rational points of a $K$-split, simple adjoint algebraic group defined over $K$. Then the parahoric 
subgroups of $\GG$ are defined. We define $P_\JJ\sub\GG$ $(\JJ\sneq\II)$ as in 3.1. 

It is likely that the definitions, results and conjectures in \S3 extend to the present $\GG$. One of the 
ingredients for such an extension is the definition of a structure of inductive limit of projective 
varieties over $\kk$ for $\GG/P_\em$. This is what we will try to achieve in this section, using as a model
the construction of the affine Grassmannian given in \cite{\SING, \S11}. But in the present case it seems to
be necessary to enlarge the category of projective varieties over $\kk$ by declaring that the inverse of a 
Frobenius map is a morphism; the inductive limit above will be taken in this enlarged category.
(The case where $\GG$ is of type $A$ has been previously considered in \cite{\ZH}. I thank X. Zhu 
for pointing out this reference.)

Let $\fg$ be the Lie algebra (over $K$) of $\GG$ with bracket $[,]$ and with Killing form $(,)$. Let 
$N=\dim\fg$. We shall assume that $p$ is such that the Killing form of the Lie algebra of a simple adjoint 
group over $\kk$ of the same type as $\GG$ is nondegenerate. For any $\co$-lattice $L$ in $\fg$ we set 
$L^\sha=\{x\in\fg;(x,L_0)\sub\co\}$; this is again an $\co$-lattice in $\fg$. We can find an $\co$-lattice 
$L_0$ in $\fg$ such that $[L_0,L_0]\sub L_0$, $L_0^\sha=L_0$ and $\{g\in\GG;\Ad(g)L_0=L_0\}=P_{\JJ_0}$ for 
some $\JJ_0\sub\II$ such that $|\JJ_0|=|\II|-1$.

Let $X$ be the set of all $\co$-lattices $L$ in $\fg$ such that $L^\sha=L$ and $[L,L]\sub L$. As in 
\cite{\SING, \S11} we have a bijection $\GG/P_{\JJ_0}@>\si>>X$, $gP_{\JJ_0}\m\Ad(g)L_0$. For any $n\in\NN$ 
let $X_n=\{L\in X;p^nL_0\sub L\sub p^{-n}L_0\}$.

Let $V_n=p^{-n}L_0/p^nL_0$, $V'_n=p^{-n}L_0/p^{2n}L_0$. Note that $V_n$ is a free $\co/p^{2n}\co$-module of 
rank $N$ and $V'_n$ is a free $\co/p^{3n}\co$-module of rank $N$. Since $p^{2n}L_0\sub p^nL_0$, we have an 
obvious map $\p_n:V'_n@>>>V_n$.

There is a well defined symmetric $\co$-bilinear pairing $(?,?)_n:V_n\T V_n@>>>p^{-2n}\co/\co$ induced by 
the restriction of the Killing form to $p^{-n}L_0$. There is a well defined alternating $\co$-trilinear map 
$(?,?,?)_n:V'_n\T V'_n\T V'_n@>>>p^{-3n}\co/\co$ induced by the map 
$p^{-n}L_0\T p^{-n}L_0\T p^{-n}L_0@>>>p^{-3n}\co$, $x,y,z\m([x,y],z)$.

For any torsion $\co$-module $Z$ we write $d(Z)=k_1+k_2+\do+k_r$ where $Z$ is isomorphic to
$\op_{i=1}^k\co/p^{k_i}\co$.

Let $E$ be the set of $\co$-lattices $L\sub\fg$ such that $p^nL_0\sub L\sub p^{-n}L_0$ and $L^\sha=L$. Let 
$E'$ be the set of $\co$-submodules $Z$ of $V_n$ such that $d(Z)=nN$ and $(Z,Z)_n=0$. We show

(a) {\it $L\m L/p^nL_0$ is a bijection $E@>\si>>E'$.}
\nl
Let $L\in E$. We have $d(L/p^nL_0)=d((p^nL_0)^\sha/L^\sha)=d(p^{-n}L_0/L)$. Hence
$$2d(L/p^nL_0)=d(L/p^nL_0)+d(p^{-n}L_0/L)=d(p^{-n}L_0/p^nL_0)=2nN$$
so that $d(L/p^nL_0)=nN$. It is clear that $(L/p^nL_0,L/p^nL_0)_n=0$. Thus the map in (a) is well defined.

Now let $Z\in E'$. Let $L$ be the inverse image of $Z$ under the obvious map $p^{-n}L_0@>>>V_n$. Note that 
$L$ is an $\co$-lattice in $\fg$ and we have clearly $(L,L)\sub\co$ so that $L\sub L^\sha$. Since 
$p^nL_0\sub L\sub p^{-n}L_0$ we have $(p^{-n}L_0)^\sha\sub L^\sha\sub(p^nL_0)^\sha$ hence 
$p^nL_0\sub L^\sha\sub p^{-n}L_0$. We have $L/p^nL_0=Z$ and 
$$\align&nN=d(L/p^nL_0)=d((p^nL_0)^\sha/L^\sha)=d(p^{-n}L_0/L^\sha)\\&=d(p^{-n}L_0/p^nL_0)-d(L^\sha/p^nL_0)=
2nN-d(L^\sha/p^nL_0).\endalign$$
Thus $nN=d(L^\sha/p^nL_0)=d(L^\sha/L)+d(L/p^nL_0)=d(L^\sha/L)+nN$ so that $d(L^\sha/L)=0$ and $L^\sha=L$. We
see that $L\in E$. The two maps $E@>>>E'$, $E'@>>>E$ defined above are clearly inverse to each other. Thus 
(a) holds.

Let $E'_0$ be the set of all $Z\in E'$ such that, setting $Z_1=\p_n\i(Z)$ (an
 $\co$-submodule of $V'_n$) we 
have $(Z_1,Z_1,Z_1)_n=0$. We show:

(b) {\it $L\m L/p^nL_0$ is a bijection $X_n@>\si>>E'_0$.}
\nl
Let $L\in X_n$. We set $Z=L/p^nL_0$. By (a) we have $Z\in E'$. Then $Z_1:=\p_n\i(Z)=L/p^{2n}L_0$ clearly
satisfies $(Z_1,Z_1,Z_1)_n=0$. Thus $Z\in E'_0$. Conversely, let $Z\in E'_0$. Let $L$ be the inverse image
of $Z$ under the obvious map $p^{-n}L_0@>>>V_n$. By (a) we have $L\in E$. Let $Z_1=L/p^{2n}L_0$. Let 
$x,y\in L$. Since $Z\in E'$ for any $z\in L$ we have $([x,y],z)\in\co$. Hence $[x,y]\in L^\sha$. Since 
$L^\sha=L$ it follows that $[x,y]\in L$. Thus $[L,L]\sub L$ so that $L\in X_n$. The result follows from (a).

Let $Y_n$ be the set of all $\co$-submodules $Z$ of $V_n$ such that $d(Z)=nN$. Let $Y'_n$ be the set of all 
$\co$-submodules $Z_1$ of $V'_n$ such that $d(Z_1)=2nN$. Note that $Z\m\p_n\i(Z)$ defines a map 
$\nu:Y_n@>>>Y'_n$. In \cite{\ZH}, X. Zhu shows that $Y_n$, $Y'_n$ 
have natural structures of projective varieties over $\kk$. 

Let $\tY_n=\{Z\in Y_n;(Z,Z)_n=0\}$, $\tY'_n=\{Z_1\in Y'_n;(Z_1,Z_1,Z_1)_n=0\}$. From the definitions, 
$\tY_n$ is a closed subvariety of $Y_n$, $\tY'_n$ is a closed subvariety of $Y'_n$ and the map 
$\nu:Y_n@>>>Y'_n$, $Z\m\p_n\i(Z)$ is a morphism of algebraic varieties. Hence 
$X'_n=\nu\i(\tY'_n)\cap\tY_n$ is a closed subvariety of $Y_n$. We can reformulate (b) as follows:

(c) $L\m L/p^nL_0$ is a bijection $X_n@>\si>>X'_n$.
\nl
It follows that $X_n$ has a natural structure of projective variety over $\kk$. Note that 
$X_0\sub X_1\sub X_2\sub\do$ and the inclusions are imbeddings of projective varieties (in the enlarged 
category above). Hence $X=\cup_{n\ge0}X_n$ is naturally an inductive limit of projective varieties over 
$\kk$.

\subhead 5.2\endsubhead
For each $n\in\NN$, $X_n$ carries a vector bundle whose fibre at $L$ is $L/pL$ viewed as a Lie algebra over 
$\kk$ in which there is a well defined notion of Borel subalgebra. We can then form the set $\tX_n$ of pairs
$(L,\fb)$ where $L\in X_n$ and $\fb$ is a Borel subalgebra of $L/pL$. This is a projective variety over 
$\kk$ which fibres over $X_n$ and each fibre is an (ordinary) flag manifold of some $L/pL$. We have 
naturally $\tX_0\sub\tX_1\sub\tX_2\sub\do$ and the inclusions are imbeddings of projective varieties (in the
enlarged category above). Hence $\tX:=\cup_{n\ge0}\tX_n$ is naturally an inductive limit of projective
varieties over $\kk$. Note that $\tX$ is the set of pairs $(L,\fb)$ where $L\in X$ and $\fb$ is a Borel 
subalgebra of $L/pL$. Let $\fb_0$ be a 
Borel subalgebra of $L_0/pL_0$. Then $gP_\em\m(\Ad(g)L_0,\Ad(g)\fb_0)$ 
is a bijection $\GG/P_\em@>\si>>\tX$. We see that $\GG/P_\em$ is naturally an inductive limit of projective 
varieties over $\kk$.


\widestnumber\key{KmL1}
\Refs
\ref\key\KL\by D.Kazhdan and G.Lusztig\paper Fixed point varieties on affine flag manifolds\jour Isr.J.Math.
\vol62\yr1988\pages129-168\endref
\ref\key\Kml\by J.L.Kim and G.Lusztig\paper On the Steinberg character of a semisimple p-adic group\lb\jour
arxiv:1204.4712\endref
\ref\key\KmL\by J.L.Kim and G.Lusztig\paper On the character of unipotent representations of a semisimple 
p-adic group\jour arxiv:1208.0320\endref
\ref\key\FOUR\by G.Lusztig\paper Unipotent representations of a finite Chevalley group of type $E_8$\jour
Quart.J.Math.\vol30\yr1979\pages315-338\endref
\ref\key\SING\by G.Lusztig\paper Singularities, character formulas and a $q$-analog of weight multiplicities
\jour Ast\'erisque\vol101-102\yr1983\pages208-229\endref
\ref\key\ICC\by G.Lusztig\paper Intersection cohomology complexes on a reductive group\jour Invent.Math.
\vol75\yr1984\pages205-272\endref
\ref\key\CUS\by G.Lusztig\paper Cuspidal local systems and graded Hecke algebras I\jour Publ.Math. I.H.E.S.
\vol67\yr1988\pages145-202\endref
\ref\key\CUSII\by G.Lusztig\paper Cuspidal local systems and graded Hecke algebras II\inbook
Representations of groups\bookinfo ed. B.Allison, Canad.Math.Soc.Conf.Proc.\vol16\publ Amer.Math.Soc. 
\year1995\pages217-275\endref
\ref\key\CLA\by G.Lusztig\paper Classification of unipotent representations of simple $p$-adic groups\lb\jour
Int.Math.Res.Notices\yr1995\vol517-589\endref
\ref\key\CLAII\by G.Lusztig\paper Classification of unipotent representations of simple $p$-adic groups II
\jour Represent.Th.\vol6\yr2002\pages243-289\endref
\ref\key\AFF\by G.Lusztig\paper Affine Weyl groups and conjugacy classes in Weyl groups\jour Transfor.
Groups\vol1\yr1996\pages83-97\endref
\ref\key\CDG\by G.Lusztig\paper Character sheaves on disconnected groups, I\jour Represent.Th.\vol7\yr2003
\pages374-403\moreref II\vol8\yr2004\pages72-124\moreref III\vol8\yr2004\pages125-144\moreref IV\vol8\yr2004
\pages145-178\moreref V\vol8\yr2004\pages346-376\moreref VI\vol8\yr2004\pages377-413\moreref VII\vol9\yr2005
\pages209-266\moreref VIII\vol10\yr2006\pages314-352\moreref IX\vol10\yr2006\pages353-379\moreref X\vol13
\yr2009\pages82-140\endref
\ref\key\CLE\by G.Lusztig\paper On the cleanness of cuspidal character sheaves\jour Moscow Math.J.\vol12
\yr2012\pages621-631\endref
\ref\key\YU\by Z.Yun\paper The spherical part of the local and global Springer actions\jour arxiv:1106.2259
\endref
\ref\key\ZH\by X.Zhu\paper Affine Grassmannians and the geometric Satake in mixed characteristic\jour
arxiv:1407.8519\endref
\endRefs
\enddocument